\documentclass[a4paper,12pt]{article}

\usepackage{float}
\usepackage{amsmath} 
\usepackage{amssymb}
\usepackage{amsfonts}
\usepackage{amsthm}
\usepackage{epstopdf}
\usepackage{color,fancybox,graphicx}
\usepackage{url}
\usepackage{psfrag}
\usepackage{color}
\usepackage{pifont}
\usepackage{mathrsfs}
\usepackage{dsfont}
\usepackage{algorithm, algorithmic}
\usepackage{enumerate}

\usepackage{pstricks, pst-tree}

\numberwithin{equation}{section}

\newtheorem{lem}{Lemma}[section]

\newtheorem{pro}{Proposition}[section]
\newtheorem{theo}{Theorem}[section]

\newcommand{\bX}{\mathbf{X}}

\newcommand{\btheta}{\boldsymbol{\hat \theta}}

\definecolor{gris25}{gray}{0.90}

\usepackage{natbib}
\bibpunct{(}{)}{;}{a}{,}{,}

\usepackage{hyperref}
\hypersetup{
colorlinks = true,
urlcolor = blue, 
linkcolor = blue,
citecolor = blue,
}
\usepackage{enumerate}

\begin{document}

\begin{center}
{\Large 
\textbf{\textsf{The Statistical Performance of Collaborative Inference}}}
\medskip
\medskip
\end{center}

{\bf G\'erard Biau}\\
{\it Sorbonne Universit\'es, UPMC Univ Paris 06, F-75005, Paris, France\\
\& Institut universitaire de France}\\
\href{mailto:gerard.biau@upmc.fr}{gerard.biau@upmc.fr}
\bigskip

{\bf Kevin Bleakley}\\
{\it INRIA Saclay, France\\
\& D\'epartement de Math\'ematiques d'Orsay, France}\\
\href{mailto:kevin.bleakley@inria.fr}{kevin.bleakley@inria.fr}
\bigskip

{\bf Beno\^{\i}t Cadre}\\
{\it IRMAR, ENS Rennes, CNRS, UEB, France}\\
\href{mailto:benoit.cadre@ens-rennes}{benoit.cadre@ens-rennes.fr}
\bigskip

\begin{abstract}
\noindent {\rm The statistical analysis of massive and complex data sets will require the development of algorithms that depend on distributed computing and collaborative inference. Inspired by this, we propose a collaborative framework that aims to estimate the unknown mean $\theta$ of a random variable $X$. In the model we present, a certain number of calculation units, distributed across a communication network represented by a graph, participate in the estimation of $\theta$ by  sequentially receiving independent data from $X$ while exchanging messages via a stochastic matrix $A$ defined over the graph. 
We give precise conditions on the matrix $A$ under which the statistical precision of the individual units is comparable to that of a (gold standard) virtual centralized estimate, even though each unit does not have access to all of the data. We show in particular the fundamental role played by both the non-trivial eigenvalues of $A$ and the Ramanujan class of expander graphs, which provide remarkable performance for moderate algorithmic cost. 
\medskip
 
\noindent \emph{Index Terms} --- Distributed computing, collaborative estimation, sto\-ch\-as\-tic matrix, graph theory, complexity, Ramanujan graph.
\medskip

\noindent \emph{2010 Mathematics Subject Classification}: 62F12, 68W15.}

\end{abstract}

\section{Introduction}
\setcounter{MaxMatrixCols}{20}
A promising way to overcome computational problems associated with inference and prediction in large-scale settings is to take advantage of distributed and collaborative algorithms, whereby several processors perform computations and exchange messages with the end-goal of minimizing a certain cost function. For instance, in modern data analysis one is frequently faced with problems where the sample size is too large for a single computer or standard computing resources. Distributed processing of such large data sets is often regarded as a possible solution to data overload, although designing and analyzing algorithms in this setting is challenging. Indeed, good distributed and collaborative architectures should maintain the desired statistical accuracy of their centralized counterpart, while retaining sufficient flexibility and avoiding communication bottlenecks which may excessively slow down  computations. The literature is too vast to permit anything like a fair summary within the confines of a short introduction---the papers by  \cite{DuAgWa12}, \cite{Jo13}, \citet{ZhDuWa13}, and references therein contain a sample of relevant work. 

Similarly, the advent of sensor, wireless and peer-to-peer networks in science and technology necessitates the design of distributed and information-exchange algorithms \citep{BoGhPrSh06,PrKuPo09}. Such networks are designed to perform inference and prediction tasks for the environments they are sensing. Nonetheless, they are typically characterized by constraints on energy, bandwidth and/or privacy, which limit the sensors' ability to share data with each other or with a hub for centralized processing. For example, in a hospital network, the aim is to make safer decisions by sharing information between therapeutic services. However, a simple exchange of database entries containing patient details can pose information privacy risks. At the same time, a large percentage of medical data may require exchanging high-resolution images, the centralized processing of which may be computationally prohibitive. Overall, such constraints call for the design of communication-constrained distributed procedures, where each node exchanges information with only a few of its neighbors at each time instance. The goal in this setting is to distribute the learning task in a computationally efficient way, and make sure that the statistical performance of the network matches that of the centralized version.

The foregoing observations have motivated the development and analysis of many local message-passing algorithms for distributed and collaborative inference, optimization and learning. Roughly speaking, message-passing procedures are those that use only local communication to approximately achieve the same end as global (i.e., centralized) algorithms, which require sending raw data to a central processing facility. Message-passing algorithms are thought to be efficient by virtue of their exploitation of local communication. They have been successfully involved in kernel linear least-squares regression estimation \citep[][]{PrKuPo09}, support vector machines \citep[][]{FoCaGi10}, sparse $L_1$ regression \citep[][]{MaBaGi10}, gradient-type optimization \citep[][]{TsBeAt86,BeTs89}, and various online inference and learning tasks \citep{BiFoHaJa11,BiFoHaJa11-2,BiClJaMo13}. An important research effort has also been devoted to so-called averaging and consensus problems, where a set of autonomous agents---which may be sensors or nodes of a computer network---compute the average of their opinions in the presence of restricted communication capabilities and try to agree on a collective decision \citep[e.g.,][]{BlHeOlTs05,OlTs11}. 

However, despite their rising success and impact in machine learning, little is known regarding the statistical properties of message-passing algorithms. The statistical performance of collaborative computing has so far been studied in terms of consensus (i.e., whether all nodes give the same result), with perhaps mean convergence rates \citep[e.g.,][]{OlTs11,DuAgWa12,ZhDuWa13}. While it is therefore proved that using a network, even sparse (i.e., with few connections), does not degrade the rate of convergence, the problem of whether it is optimal to do this remains unanswered, including for the most basic statistics. For example, which network properties guarantee collaborative calculation performances equal to those of a hypothetical centralized system? The goal of this article is to give a more precise answer to this fundamental question. In order to present in the clearest way possible the properties such a network must have, we undertake this study for the most simple statistic possible: the mean.

In the model we consider, there are a number of computing agents (also known as nodes or processors) that sequentially estimate the mean of a random variable by regularly updating an estimate stored in their memory. Meanwhile, they exchange messages, thus informing each other about the results of their latest computations. Agents that receive messages use them to directly update  the value in their memory by forming a convex combination. We focus primarily on the properties that the communication process must satisfy to ensure that the statistical precision of a single processor---that only sees part of the data---is similar to that of an inaccessible centralized intelligence that could tackle the whole data set at once. The literature is surprisingly quiet on this question, which we believe is of fundamental importance if we want to provide concrete tradeoffs between communication constraints and statistical accuracy. 

This paper makes several important contributions. First, in Section 2 we introduce communication network models and define a performance ratio allowing us to quantify the statistical quality of a network. In Section 3 we  analyze the asymptotic behavior of this performance ratio as the number of data items $t$ received online sequentially  per node becomes large, and give precise conditions on communication matrices $A$ so that this ratio is asymptotically optimal.
Section 4 goes one step further, connecting the rate of convergence of the ratio with the behavior of the eigenvalues of $A$. In Section 5 we present the remarkable Ramanujan expander graphs and analyze the tradeoff between statistical efficiency and communication complexity for these graphs with a series of simulation studies. Lastly, Section 6 provides several elements for  analysis of more complicated asynchronous models with delays. For clarity, proofs are gathered  in Section 7.
\section{The model}
Let $X$ be a square-integrable real-valued random variable, with $\mathbb E X=\theta$ and $\mbox{Var} (X)=\sigma^2$. We consider a set $\{1, \hdots, N\}$ of computing entities $(N\geq 2)$ that collectively participate in the estimation of $\theta$. In this distributed model, agent $i$ sequentially receives an i.i.d.~sequence $X^{(i)}_1, \hdots, X^{(i)}_t, \hdots,$ distributed as the prototype $X$, and forms, at each time $t$, an estimate of $\theta$. It is assumed throughout that the $X_t^{(i)}$ are independent when both $t \geq 1$ and $i\in \{1, \hdots, N\}$ vary. 

In the absence of communication between agents, the natural estimate held by agent $i$ at time $t$ is the empirical mean 
$$\bar X_{t}^{(i)}=\frac{1}{t}\sum_{k=1}^t X_k^{(i)}.$$
Equivalently, processor $i$ is initialized with $X^{(i)}_1$ and performs its estimation via the iteration
$$\bar X^{(i)}_{t+1}=\frac{t \bar X_t^{(i)}+X_{t+1}^{(i)}}{t+1}, \quad t\geq 1.$$
Let $\top$ denote transposition and assume that vectors are in column format. Letting $\bX_t=(X_t^{(1)}, \hdots, X_t^{(N)})^\top$ and $\bar \bX_t=(\bar X_t^{(1)}, \hdots, \bar X_t^{(N)})^\top$, we see that 
\begin{equation}
\label{nocom}
\bar \bX_{t+1}=\frac{t \bar \bX_t+\bX_{t+1}}{t+1},\quad t\geq 1.
\end{equation}
In a more complicated collaborative setting, besides its own measurements and computations, each agent may also receive messages from other processors and combine this information with its own conclusions. At its core, this message-passing process can be modeled by a {\it directed} graph $\mathscr G=(\mathscr V,\mathscr E)$ with vertex set $\mathscr V=\{1,\hdots,N\}$ and edge set $\mathscr E$. This graph represents the way agents communicate, with an edge from $j$ to $i$ (in that order) if $j$ sends  information to $i$. Furthermore, we have an $N\times N$ stochastic matrix $A=(a_{ij})_{1\leq i,j\leq N}$  (i.e., $a_{ij}\geq 0$ and for each $i$, $\sum_{j=1}^Na_{ij}=1$)
with associated graph $\mathscr G$, i.e., $a_{ij}>0$ if and only if $(j,i)\in\mathscr E$. The matrix $A$ accounts for the way agents incorporate  information during the collaborative process. Denoting by $\btheta_t=(\hat \theta^{(1)}_t, \hdots, \hat \theta^{(N)}_t)^\top$ the collection of estimates held by the $N$ agents over time, the computation/combining mechanism is assumed to be as follows:
$$\btheta_{t+1}=\frac{t}{t+1}A\btheta_t+\frac{1}{t+1}\bX_{t+1}, \quad t\geq 1,$$
with $\btheta_{1}=(X_1^{(1)}, \hdots, X_1^{(N)})^\top$. Thus, each individual estimate $\hat \theta^{(i)}_{t+1}$ is a convex combination of the estimates $\hat \theta^{(j)}_t$ held by the agents over the network at time $t$, augmented by the new observation $X^{(i)}_{t+1}$. 

The matrix $A$ models the way processors exchange messages and collaborate, ranging from $A=I_N$ (the $N\times N$ identity matrix, i.e., no communication) to $A=\mathbf 1 \mathbf 1^{\top} /N$ (where $\mathbf 1=(1, \hdots, 1)^{\top}$, i.e., full communication). We note in particular that the choice $A=I_N$ gives back iteration (\ref{nocom}) with $\btheta_t=\bar \bX_t$. We also note that, given a graph $\mathscr G$, various choices are possible for $A$. Thus, aside from a convenient way to represent a communication channel over which agents can retrieve information from each other, the matrix $A$ can be seen as a ``tuning parameter'' on $\mathscr G$ to improve the statistical performance of $\btheta_t$, as we shall see later. Important examples for $A$ include the choices
\begin{equation}
\label{A1}
A_1=\frac{1}{2}
\begin{pmatrix}
   1 & 1 & & &\\
    1 & 0 & 1 & &\\
    & 1 &0 & 1 &\\
 \vdots &  \vdots & \vdots & \vdots & \vdots &\vdots& \vdots & \vdots &\vdots\\
  & & & & &1 & 0 &1 &\\
 & & & & & & 1 &0 &1\\
 & & & & & & & 1&1 
\end{pmatrix}
\end{equation}
and
\begin{equation}
\label{A2}
A_2=\frac{1}{3}
\begin{pmatrix}
   2 & 1 &  & \\
   1 & 1 & 1 &  &\\
    & 1 & 1 & 1& & \\
       \vdots &\vdots& \vdots & \vdots &\vdots&\vdots&\vdots&\vdots &\vdots\\
   &&  & &  &1 &1&1& \\
      &&  & &  & &1&1& 1\\
       &&  & &  & &&1& 2\\
        \end{pmatrix}
\end{equation}
(unmarked entries are zero). It is easy to verify that for all $t\geq 1$, 
\begin{equation}
\label{BLVDP}
\btheta_{t}=\frac{1}{t}\sum_{k=0}^{t-1} A^k \bX_{t-k}.
\end{equation}
Thus, denoting by $\|\cdot\|$ the Euclidean norm (for vector or matrices), we may write, for all $t\geq 1$,
\begin{align*}
\mathbb E \|\btheta_{t}-\theta \mathbf 1\|^2 &=\frac{1}{t^2}\mathbb E \bigg\| \sum_{k=0}^{t-1}A^k(\bX_{t-k}-\theta \mathbf 1)\bigg\|^2\\
& \quad \mbox{(since $A^k$ is a stochastic matrix)}\\
& =\frac{1}{t^2}\sum_{k=1}^{t}\mathbb E \left\| A^{t-k}(\bX_{k}-\theta \mathbf 1)\right\|^2,
\end{align*}
by independence of $\bX_1, \hdots, \bX_t$. It follows that
\begin{align*}
\mathbb E \|\btheta_{t}-\theta \mathbf 1\|^2  &\leq \mathbb E\|\bX_{1}-\theta\mathbf 1\|^2 \times \frac{1}{t^2}\sum_{k=0}^{t-1} \| A^{k}\|^2\\
& \leq \mathbb E\|\bX_{1}-\theta \mathbf 1\|^2 \times \frac{N}{t}.
\end{align*}
In the last inequality, we used the fact that $A^k$ is a stochastic matrix and thus $\|A^k\|^2\leq N$ for all $k\geq 0$. We can merely conclude that $\mathbb E \|\btheta_{t}-\theta \mathbf 1\|^2 \to 0$ as $t\to \infty$ (mean-squared error consistency), and so $\hat \theta_t^{(i)} \to \theta$ in probability for each $i\in \{1, \hdots, N\}$. Put differently, the agents asymptotically agree on the (true) value of the parameter, independently of the choice of the (stochastic) matrix $A$---this property is often called consensus in the distributed optimization literature \citep[see, e.g.,][]{BeTs89}.

The consensus property, although interesting, does not say anything about the positive (or negative) impact of the graph on the comparative performances of estimates with respect to a centralized version. To clarify this remark, assume that there exists a centralized intelligence that could tackle all data $X_1^{(1)}, \hdots, X_t^{(1)}, \hdots, X_1^{(N)}, \hdots, X_t^{(N)}$ at time $t$, and take advantage of these samples to assess the value of the parameter $\theta$. In this ideal framework, the natural estimate of $\theta$ is the global empirical mean
$$\bar{\mathbb X}_{Nt}=\frac{1}{Nt}\sum_{i=1}^N\sum_{k=1}^tX_{k}^{(i)},$$
which is clearly the best we can hope for with the data at hand. However, this estimate is to be considered as an unattainable ``gold standard'' (or oracle), insofar as it uses the whole $(N\times t)$-sample. In other words, its evaluation requires sending all examples to a centralized processing facility, which is precisely what we want to avoid. 

\begin{figure}[!h]
\begin{center}
\includegraphics[width=10cm, height = 8cm]{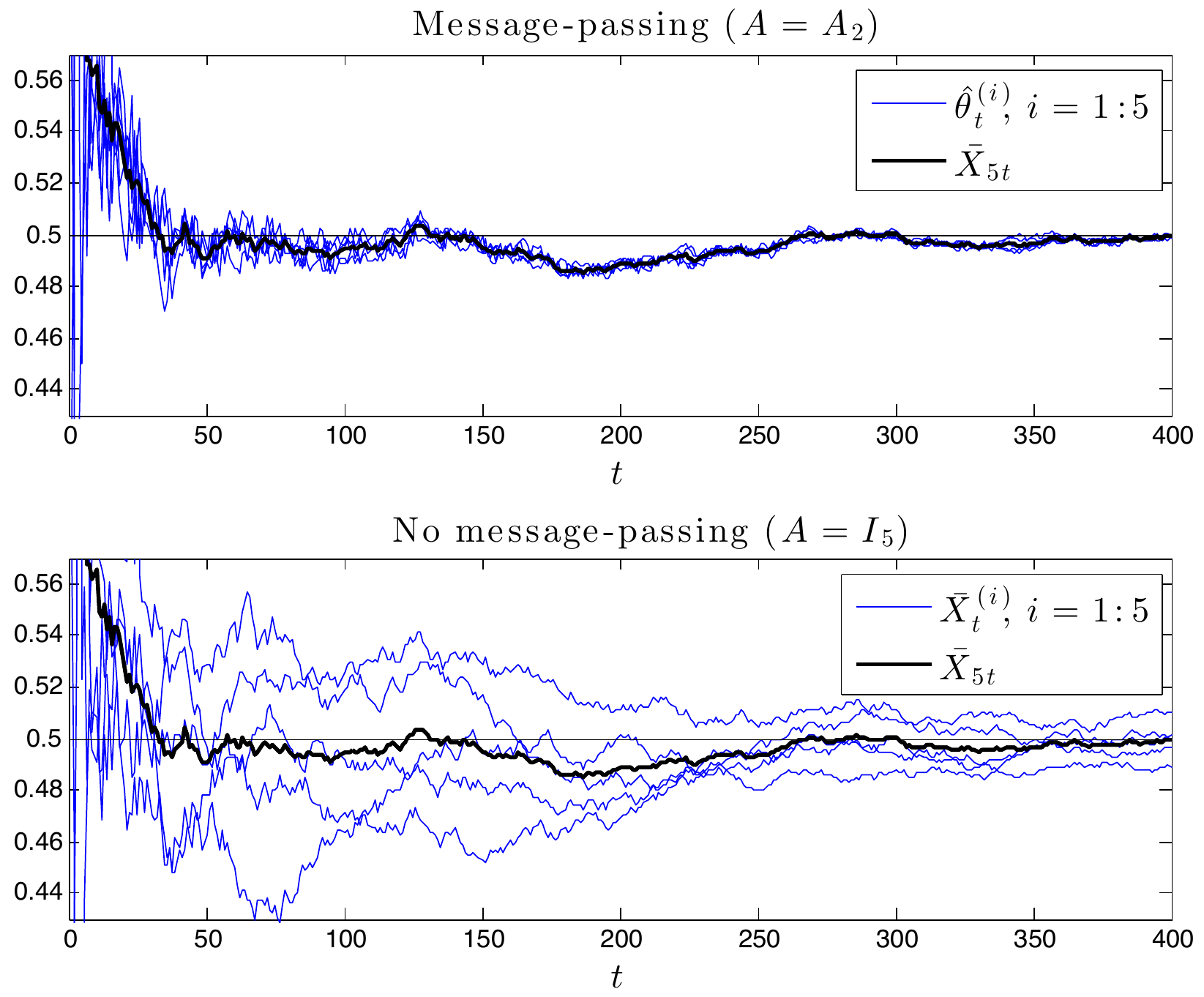}
\end{center}
\caption{Convergence of individual nodes' estimates with and without message-passing.}\label{fig5}
\end{figure}

Thus, a natural question arises: can the message-passing process be tapped to ensure that the individual estimates $\hat \theta_t^{(i)}$ achieve statistical accuracy ``close'' to that of the gold standard $\bar{\mathbb X}_{Nt}$? 
Figure~\ref{fig5}
illustrates this pertinent question.

In the  trials shown, i.i.d.~uniform random variables on $[0,1]$ are delivered online to 
$N=5$ nodes, one to each at each time $t$. With message-passing (here, $A = A_2$), each node aggregates the new data point with data it has seen previously and messages received from its nearest neighbors in the network. We see that all of the five nodes' updates seem to converge with a performance comparable to that of the (unseen) global estimate  $\bar{\mathbb X}_{Nt}$ to the mean 0.5. In contrast, in the absence of message-passing ($A = I_5$), individual nodes' estimates do still  converge to 0.5, but at a slower rate.

To deal with this question of statistical accuracy  satisfactorily, we first need a criterion to compare the performance of $\btheta_t$ with that of $\bar{\mathbb X}_{Nt}$. Perhaps the most natural one is the following ratio, which depends upon the matrix $A$:
$$\tau_t(A)=\frac{\mathbb E \left\| (\bar{\mathbb X}_{Nt}-\theta)\mathbf 1\right\|^2}{\mathbb E\|\btheta_t-\theta \mathbf 1\|^2}, \quad t \geq 1.$$
The more this ratio is close to 1, the more the collaborative algorithm is statistically efficient, in the sense that its performance compares favorably to that of the centralized gold standard. In the remainder of the paper, we call $\tau_t(A)$ the {\it performance ratio} at time $t$. 

Of particular interest in our approach is the stochastic matrix $A$, which plays a crucial role in the analysis. Roughly, a good choice for $A$ is one for which $\tau_t(A)$ is not too far from $1$, while ensuring that communication over the network is not prohibitively expensive. Although there are several ways to measure ``complexity'' of the message-passing process, we have in mind a setting where the communication load is well-balanced between agents, in the sense that no node should play a dominant role. To formalize this idea, we define the communication-complexity index $\mathscr C(A)$ as the maximal indegree of the edges of the graph $\mathscr G$ associated with $A$, i.e., the maximal number of edges pointing to a node in $\mathscr G$ (by convention, self-loops are counted twice when $\mathscr G$ is undirected). Essentially, $A$ is communication-efficient when $\mathscr C(A)$ is small with respect to $N$ or, more generally, when $\mathscr C(A)=\mbox{O}(1)$ as $N$ becomes large.

To provide some context, $\mathscr C(A)$ measures in a certain sense the ``local'' aspect of message exchanges induced by $A$. We have in mind  node connection set-ups where $\mathscr C(A)$ is small, perhaps due to energy or bandwidth constraints in the system's architecture, or when for privacy reasons data must not be sent to a central node. 
Indeed, a large $\mathscr C(A)$ roughly means that one or several nodes play centralized roles---precisely what we are trying to avoid. 
Furthermore, the decentralized networks we are interested in can be seen as being more autonomous than high-$\mathscr C(A)$ ones, in the sense that having few network connections means less things that can potentially break, as well as improved robustness due to the fact that the loss of one node does not lead to destruction of the whole system. As examples, the matrices $A_1$ and $A_2$ defined earlier have $\mathscr C(A_1)=3$ and $\mathscr C(A_2)=4$, respectively, while the stochastic matrix $A_3$ below has $\mathscr C(A_3)=N+1$:
\begin{equation}
\label{A3}
A_3=\frac{1}{N}
\begin{pmatrix}
   1 & 1 & 1& \cdots &1 & 1 & 1\\
   1  & N-1 & & & &   \\
    1       &  & N-1& & &   \\
   \vdots & \vdots & \vdots & \vdots &\vdots& \vdots & \vdots \\
 1 & & & & & & N-1
\end{pmatrix}.
\end{equation}
Thus, from a network complexity point of view, $A_1$ and $A_2$ are preferable to $A_3$ where node 1 has the flavor of a central command center. 

Now, having defined $\tau_t(A)$ and $\mathscr C(A)$, it is natural to suspect that there will be some kind of tradeoff between implementing a low-complexity message-passing algorithm (i.e., $\mathscr C(A)$ small) and achieving good asymptotic performance (i.e., $\tau_t(A)\approx 1$ for large $t$). Our main goal in the next few sections is to probe this intuition by analyzing the asymptotic behavior of $\tau_t(A)$ as $t\to \infty$ under various assumptions on $A$. We start by proving that $\tau_t(A)\leq 1$ for all $t\geq 1$, and give precise conditions on the matrix $A$ under which $\tau_t(A) \to 1$. Thus, thanks to the benefit of inter-agent communication, the statistical accuracy of individual estimates may be asymptotically comparable to that of the gold standard, despite the fact that none of the agents in the network have access to all of the data. Indeed, as we shall see, this stunning result is possible even for low-$\mathscr C(A)$ matrices. The take-home message here is that the communication process, once cleverly designed, may ``boost'' the individual estimates, even in the presence of severe communication constraints. We also provide an asymptotic development of $\tau_t(A)$, which offers valuable information on the optimal way to design the communication network in terms of the eigenvalues of $A$.
\section{Convergence of the performance ratio}
\label{cv}
Recall that a stochastic square matrix $A=(a_{ij})_{1\leq i,j\leq N}$ is irreducible if for every pair of indices $i$ and $j$, there exists a nonnegative integer $k$ such that $(A^{k})_{ij}$ is not equal to 0.  The matrix is said to be reducible if it is not irreducible.
\begin{pro}
\label{PROPO1}
We have $\frac{1}{N}\leq \tau_t(A)\leq 1$ for all $t\geq 1$. In addition, if $A$ is reducible, then
$$\tau_t(A)\leq 1-\frac{1}{N+1}, \quad t \geq 1.$$
\end{pro}

It is apparent from the proof of the proposition (all proofs are found in Section \ref{proofs}) that the lower bound $1/N$ for $\tau_t(A)$ is achieved by taking $A=I_N$, which is clearly the worst choice in terms of communication. This proposition also shows that the irreducibility of $A$ is a necessary condition for the collaborative algorithm to be statistically efficient, for otherwise there exists $\varepsilon \in (0,1)$ such that $\tau_t(A)\leq 1-\varepsilon$ for all $t\geq 1$. 

We recall from the theory of Markov chains \citep[e.g.,][]{GrSt01} that for a fixed agent $i\in\{1,\hdots,N\}$, the period of $i$ is the greatest common divisor of all positive integers $k$ such that $(A^k)_{ii}> 0$. When $A$ is irreducible, the period of every state is the same and is called the period of A. The following lemma describes the asymptotic behavior of $\tau_t(A)$ as $t$ tends to infinity.
\begin{lem}
\label{LEM1}
Assume that $A$ is irreducible, and let $d$ be its period. Then there exist projectors $Q_1,\hdots, Q_d$ such that
$$\tau_t(A) \to \frac{1}{\sum_{\ell=1}^d \|Q_{\ell}\|^2}\quad \mbox{as } t\to \infty.$$
\end{lem}

The projectors $Q_1, \hdots, Q_d$ in Lemma \ref{LEM1} originate from the decomposition
$$A^k=\sum_{\ell=1}^d \lambda_{\ell}^k Q_{\ell}+\sum_{\gamma\in\Gamma} \gamma^k Q_\gamma(k),$$
where $\lambda_1=1,\hdots,\lambda_d$ are the (distinct) eigenvalues of $A$ of unit modulus, $\Gamma$  the set of eigenvalues of $A$ of modulus strictly smaller than 1, and $Q_\gamma(k)$ certain $N\times N$ matrices (see Theorem \ref{FF} in the proofs section). In particular, we see that $\tau_t(A)\to 1$ as $t\to \infty$ if and only if $\sum_{\ell=1}^d \|Q_{\ell}\|^2=1$. It turns out that this condition is satisfied if and only if $A$ is irreducible, aperiodic (i.e., $d=1$), and bistochastic, i.e., $\sum_{i=1}^Na_{ij}=\sum_{j=1}^Na_{ij}=1$  for all $(i,j)\in \{1,\hdots,N\}^2$. This important result is encapsulated in the next theorem. 
\begin{theo}
\label{convergence}
We have $\tau_t(A)\to 1$ as $t\to \infty$ if and only if $A$ is irreducible, aperiodic, and bistochastic.
\end{theo}

Theorem \ref{convergence} offers necessary and sufficient conditions 
for the communication matrix $A$ to be asymptotically statistically efficient. Put differently, under 
the conditions of the theorem, the message-passing process conveys sufficient information to local
computations to make  individual estimates as accurate as the gold standard for large $t$. In the context of multi-agent coordination,
an example of such a communication
network is the so-called (time-invariant) equal neighbor model \citep[][]{TsBeAt86,OlTs11}, in which
\begin{equation*}
a_{ij}= \left \{ \begin{array}{ll}
1/{|N^{(i)}|}&  \mbox{if $j \in N^{(i)}$}\\
0 & \mbox{otherwise},
\end{array}
\right .
\end{equation*}
where 
$$N^{(i)}=\big\{ j \in \{1, \hdots, N\}:a_{ij} > 0\big\}$$
is the set of agents whose value is taken into account by $i$, and $|N^{(i)}|$  its cardinality.
Clearly, the communication matrix $A$ is stochastic, and also bistochastic as soon
as $A$ is symmetric (bidirectional model). Assuming in addition that the directed
graph $\mathscr G$ associated with $A$ is strongly connected means that $A$ is irreducible.
Moreover, if $a_{ii}>0$ for some $i\in \{1, \hdots, N\}$, then $A$ is also aperiodic,
so the conditions of Theorem \ref{convergence} are fulfilled. 

It is interesting to note that there exist low-$\mathscr C(A)$ matrices that meet the requirements
of Theorem \ref{convergence}. This is for instance the case of  matrices $A_1$ and $A_2$ in
(\ref{A1}) and (\ref{A2}), which are irreducible, aperiodic and bistochastic, and satisfy $\mathscr C(A)\leq 4$. Also note
that the matrix $A_3$ in (\ref{A3}), though irreducible, aperiodic and bistochastic, should be avoided because $\mathscr C(A_3)=N+1$. 

We stress that the irreducibility and aperiodicity conditions are inherent properties of the graph $\mathscr G$, not $A$, insofar as these conditions do not depend upon the actual values of the nonzero entries of $A$. This is different for the bistochasticity condition, which requires knowledge of the coefficients of $A$. In fact, as observed by \citet[][]{SiKn67}, it is not always possible to associate such a bistochastic matrix with a given directed graph $\mathscr G$. To be more precise, consider $G=(g_{ij})_{1\le i,j\le N}$, the transpose of the adjacency matrix of the graph $\mathscr G$---that is, $g_{ij}\in \{0,1\}$ and $g_{ij}=1 \Leftrightarrow (j,i)\in\mathscr E$. Then $G$ is said to have total support if, for every positive element $g_{ij}$, there exists a permutation $\sigma$ of $\{1,\hdots,N\}$ such that $j=\sigma(i)$ and $\prod_{k=1}^N g_{k\sigma(k)}>0$. The main theorem of \citet[][]{SiKn67} asserts that there exists a bistochastic matrix $A$ of the form $A=D_1GD_2$, where $D_1$ and $D_2$ are $N\times N$ diagonal matrices with positive diagonals, if and only if $G$ has total support. The algorithm to induce $A$ from $G$ is called the Sinkhorn-Knopp algorithm. It does this by generating a sequence of matrices whose rows and columns are normalized alternately. It is known that the convergence of the algorithm is linear and upper bounds have been given for its rate of convergence \citep[e.g.,][]{Kn08}. 

Nevertheless, if for some reason we face a situation where it is impossible to associate a bistochastic matrix with the graph $\mathscr G$, Proposition \ref{CORO1} below shows that it is still possible to obtain information about the performance ratio, provided $A$ is irreducible and aperiodic. 
\begin{pro}
\label{CORO1}
Assume that $A$ is irreducible and aperiodic. Then
$$\tau_t(A)\to \frac{1}{N\|\boldsymbol{\mu}\|^2} \quad \mbox{as }t \to \infty,$$
where $\boldsymbol{\mu}$ is the stationary distribution of $A$.
\end{pro}

To illustrate this result, take $N=2$ and consider the graph $\mathscr G$ with (symmetric) adjacency matrix $\mathbf 1 \mathbf 1^{\top}$ (i.e., full communication). Various stochastic matrices may be associated with $\mathscr G$, each  with a certain statistical performance. For $\alpha >1$ a given parameter, we may choose for example
\begin{equation*} 
H_{\alpha}=\frac{1}{\alpha}\begin{pmatrix}
1 & \alpha-1 \\
1 & \alpha-1 \\
\end{pmatrix}.
\end{equation*}
When $\alpha=2$, we have $\tau_t(H_2)\to 1$ by Theorem \ref{convergence}. More generally, using Proposition \ref{CORO1}, it is an easy exercise to prove that, as $t\to\infty$, 
$$\tau_t(H_{\alpha})\to \frac{\alpha^2}{2+2(\alpha-1)^2}.$$
We see that the statistical performance of the local estimates deteriorates as $\alpha$ becomes large, for in this case $\tau_t(H_{\alpha})$ gets closer and closer to $1/2$. This toy model exemplifies the role the stochastic matrix is playing as a ``tuning parameter'' to improve the performance of the distributed estimate. 
\section{Convergence rates}

Theorem \ref{convergence} gives precise conditions ensuring $\tau_t(A)=1+\mbox{o}(1)$, but  does not say anything about the rate (i.e., the behavior of the second-order term) at which this convergence occurs. It turns out that a much more informative limit may be obtained at the price of the mild additional assumption that the stochastic matrix $A$ is symmetric (and hence bistochastic).
\begin{theo}
\label{theovitesse}
Assume that $A$ is irreducible, aperiodic, and symmetric. Let $1>\gamma_2\geq \cdots \geq \gamma_N > -1$ be the eigenvalues of $A$ different from $1$. Then
$$\tau_t(A) = \frac{1}{1 + \frac{1}{t}\sum_{\ell = 2}^{N}\frac{1 - \gamma_{\ell}^{2t}}{1 - \gamma_{\ell}^2}}.$$
In addition, setting 
$$\mathscr{S}(A) = \sum_{\ell=2}^N\frac{1}{1-\gamma_{\ell}^2} \quad \mbox{and} \quad \Gamma(A)=\max_{2 \leq \ell \leq N} |\gamma_\ell |,$$
we have, for all $t\geq 1$,
$$1-\frac{\mathscr S(A)}{t}\le \tau_t(A)\le 1-\frac{\mathscr S(A)}{t}+\Gamma^{2t}(A)\frac{\mathscr S(A)}{t}+\Big(\frac{\mathscr S(A)}{t}\Big)^2.$$
\end{theo}

Clearly, we thus have 
$$t\big(1-\tau_t(A)\big)\to \mathscr S(A) \quad \mbox{as }t \to \infty.$$
The take-home message is that the smaller the coefficient $\mathscr S(A)$, the better the matrix $A$ performs from a statistical point of view. In this respect, we note that $\mathscr S(A)\ge N-1$ (uniformly over the set of stochastic, irreducible, aperiodic, and symmetric matrices). Consider the full-communication matrix
\begin{equation}
\label{A0}
A_0=\frac{1}{N}
\mathbf 1 \mathbf 1^\top,
\end{equation}
which models a saturated communication network
in which each agent shares its information with all  others. The associated communication topology, which has $\mathscr C(A_0)=N+1$, is roughly equivalent to a centralized algorithm and, as such, is considered inefficient from a computational point of view. On the other hand, intuitively, the amount of statistical information propagating through the network is large so $\mathscr{S}(A_0)$ should be small. Indeed, it is easy to see that in this case, 
$\gamma_{\ell}=0$ for all $\ell \in \{2, \hdots, N\}$ and $\mathscr S(A_0)=N-1$. Therefore, although complex in terms of communication, $A_0$ is statistically optimal.

For a comparative study of statistical performance and communication complexity of matrices, let us consider the sparser graph associated with the tridiagonal matrix $A_1$ defined in (\ref{A1}). With this choice, $\gamma_{\ell}=\cos \frac{(\ell-1) \pi}{N}$ \citep[][]{Fi72}, so that
\begin{equation}
\label{AA}
\mathscr{S}(A_1)=\sum_{\ell=1}^{N-1}\frac{1}{1-\cos^2{\frac{\ell\pi}{N}}}= \frac{N^2}{6}+\mbox{O}(N) \quad \mbox{as $N\to \infty$}.
\end{equation}
Thus, we lose a power of $N$ but now have  lower communication complexity $\mathscr C(A_1)=3$.

Let us now consider the tridiagonal matrix $A_2$ defined in (\ref{A2}). Noticing that $3A_2=2A_1+I_N$, we deduce that for the matrix $A_2$,
$\gamma_{\ell}=\frac{1}{3}+\frac{2}{3} \cos \frac{(\ell-1)\pi}{N}$, $2\leq \ell \leq N$. Thus, as $N\to\infty$,
\begin{equation}
\label{BB}
\mathscr{S}(A_2)=\frac{N^2}{9}+\mbox{O}(N).
\end{equation}
By comparing (\ref{AA}) and (\ref{BB}), we can conclude that the matrices $A_1$ and $A_2$, which are both low-$\mathscr C(A)$, are also nearly equivalent from a statistical efficiency point of view. $A_2$ is nevertheless preferable to $A_1$, which has a larger constant in front of the $N^2$. This slight difference may be due to the fact that most of the diagonal elements of $A_1$ are zero, so that agents $i \in \{2,\hdots,N-1\}$ do not integrate their current value in the next iteration, as happens for $A_2$. Furthermore, for large $N$, the performance of $A_1$ and $A_2$ are expected to dramatically deteriorate in comparison with those of $A_0$, since $\mathscr S(A_1)$ and $\mathscr S(A_2)$ are proportional to $N^2$, while $\mathscr S(A_0)$ is proportional to $N$.

\begin{figure}[h]
\begin{center}
\includegraphics[width=10cm, height = 8cm]{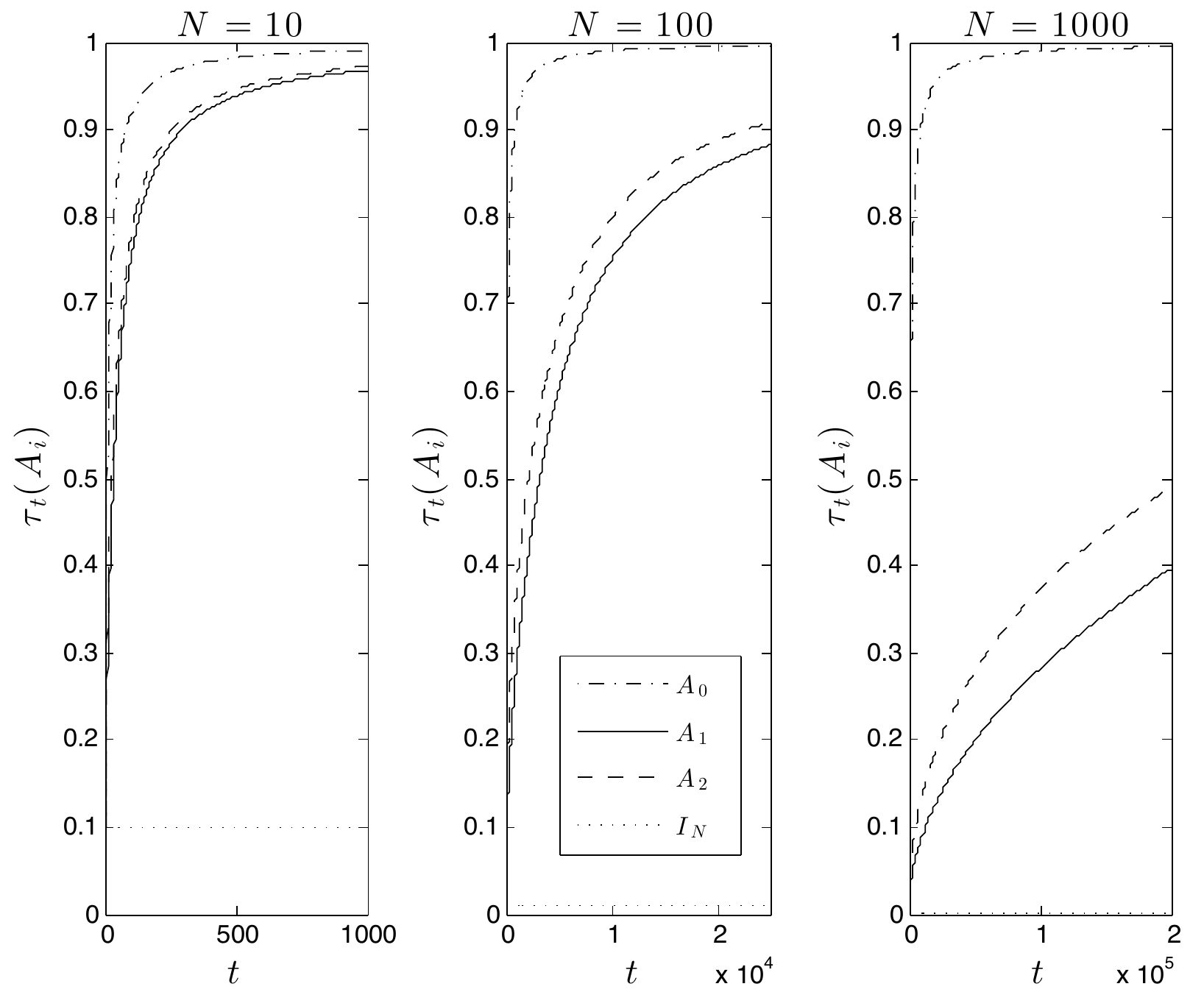}
\end{center}
\caption{Evolution of $\tau_t(A_i)$ with $t$ for different values of $N$, for $A = A_0$, $A_1$, $A_2$ and $I_N$.}\label{fig1}
\end{figure}

Figure~\ref{fig1} shows the evolution of $\tau_t(A)$ for $N$ fixed and $t$ increasing for the matrices $A = A_0$, $A_1$, $A_2$ as well as the identity  $I_N$.

As expected, we see  convergence of $\tau_t(A_i)$ to $1$, with degraded performance as the number of agents $N$ increases. Also, we see that the lack of message-passing for $I_N$ means it is statistically inefficient, with constant $\tau_t(I_N) = 1/N$ for all $t$.

The discussion and plots above highlight the crucial influence of $\mathscr S(A)$ on the performance of the communication network. Indeed, Theorem \ref{theovitesse} shows that the optimal order for $\mathscr S(A)$ is $N$, and that this scaling is achieved by the computationally-inefficient choice $A_0$---see (\ref{A0}). Thus, a natural question to ask is whether there exist communication networks that have $\mathscr S(A)$ proportional to $N$ and, simultaneously, $\mathscr C(A)$ constant or small with respect to $N$. These two conditions, which are in a sense contradictory, impose that the absolute values of the non-trivial eigenvalues $\gamma_\ell$ stay far from 1, while the maximal indegree of the graph $\mathscr G$ remains moderate. It turns out that these requirements are satisfied by so-called Ramanujan graphs, which are presented in the next section.

\section{Ramanujan graphs}
In this section, we consider {\it undirected} graphs $\mathscr G=(\mathscr V,\mathscr E)$ that are also $d$-regular, in the sense that all vertices have the same degree $d$; that is each vertex is incident to exactly $d$ edges. Recall that in this definition, self-loops are counted twice and multiple edges are allowed. However, in what follows, we restrict ourselves to graphs without self-loops and multiple edges. In this setting, the natural (bistochastic) communication matrix $A$ associated with $\mathscr G$ is $A=\frac{1}{d} G$, where $G=(g_{ij})_{1\le i,j\le N}$ is the adjacency matrix of $\mathscr G$ ($g_{ij}\in \{0,1\}$ and $g_{ij}=1 \Leftrightarrow (i,j)\in\mathscr E$). Note that $\mathscr C (A)=d$.

The matrix $G$ is symmetric and we let $d=\mu_1\geq \mu_2\geq \cdots \geq \mu_N \geq -d$ be its (real) eigenvalues. Similarly, we let $1=\gamma_1\geq \gamma_2\geq \cdots \geq \gamma_N \geq -1$ be the eigenvalues of $A$, with the straightforward correspondence $\gamma_i=\mu_i/d$. We note that $A$ is irreducible (or, equivalently, that $\mathscr G$ is connected) if and only if $d>\mu_2$ \citep[see, e.g.,][Section 2.3]{ShLiWi06}. In addition, $A$ is aperiodic as soon as  $\mu_N>-d$. According to the Alon-Boppana theorem \citep[][]{Ni91} one has, for every $d$-regular graph,
$$\mu_2 \geq 2\sqrt{d-1}-\mbox{o}_N(1),$$
where the $\mbox{o}_N(1)$ term is a quantity that tends to zero for every fixed $d$ as $N\to \infty$. Moreover, a $d$-regular graph $\mathscr G$ is called Ramanujan if 
$$\max\big(|\mu_{\ell}|:\mu_{\ell}<d\big)Ê\leq 2\sqrt{d-1}.$$
In view of the above, a Ramanujan graph is optimal, at least
as far as the spectral gap measure of expansion is concerned. Ramanujan graphs fall in the category of so-called expander graphs, which have the apparently contradictory features of being both highly connected and at the same time sparse \citep[for a review, see][]{ShLiWi06}. 

Although the existence of Ramanujan graphs for any degree larger than or equal to $3$ has been recently established by \citet{MaSpSr15}, their explicit construction remains difficult to use in practice. However, a conjecture by \citet{Al86}, proved by \citet{Fr08} \citep[see also][]{Bo15} asserts that most $d$-regular graphs are Ramanujan, in the sense that for every $\varepsilon >0$, 
$$\mathbb P\Big (\max \big(|\mu_2|,|\mu_N|\big) \geq 2\sqrt{d-1}+\varepsilon\Big)\to 0 \quad \mbox{as }N\to\infty,$$
or equivalently, in terms of the eigenvalues of $A$,
$$\mathbb P\Big (\max \big(|\gamma_2|,|\gamma_N|\big) \geq \frac{2\sqrt{d-1}}{d}+\varepsilon\Big)\to 0 \quad \mbox{as }N\to\infty.$$
In both results, the limit is along any sequence going to infinity with $Nd$ even, and the probability is with respect to random graphs uniformly sampled in the family of $d$-regular graphs with vertex set $\mathscr V=\{1, \hdots,N\}$.

\begin{figure}[h]
\begin{center}
\includegraphics[width=6cm]{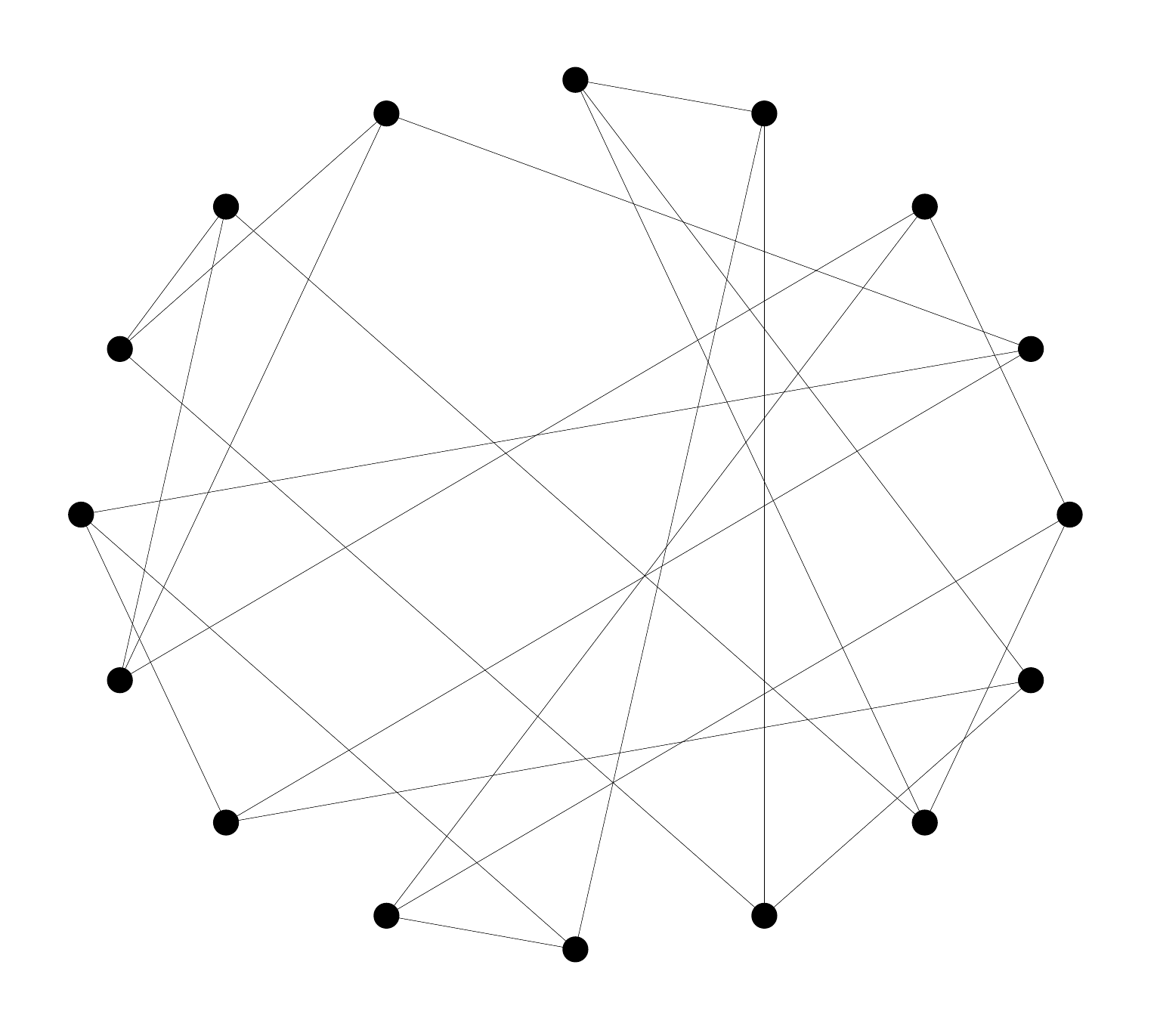}
\end{center}
\caption{Randomly-generated $3$-regular Ramanujan graph with $N=16$ vertices.}\label{ramagraph}
\end{figure}

In order to generate a random irreducible, aperiodic $d$-regular Ramanujan graph, we can first generate a random $d$-regular graph using an improved version of the standard \emph{pairing} algorithm, proposed by \citet{St99}. We retain it if it passes the tests of being irreducible, aperiodic and Ramanujan as described above. Otherwise, we
continue to generate a $d$-regular graph until all these conditions are satisfied.
Figure~\ref{ramagraph} gives an example of a $3$-regular Ramanujan graph with $N=16$ vertices, generated in this way. 

Now, given an irreducible and aperiodic communication matrix $A$ associated with a $d$-regular Ramanujan graph $\mathscr G$, we have, whenever $d\geq 3$,
$$\mathscr S (A) \leq \frac{N-1}{1-\frac{4(d-1)}{d^2}}.$$
Thus, recalling that $\mathscr S(A)\geq N-1$, we see that $\mathscr S (A)$ scales optimally as $N$ while having $\mathscr C(A)=d$ (fixed). This remarkable superefficiency property can be compared with the full-communication matrix
$A_0$, which has $\mathscr S(A_0)=N-1$ but  inadmissible complexity $\mathscr C(A_0)=N+1$.

\begin{figure}[h]
\begin{center}
\includegraphics[width=10cm, height = 8cm]{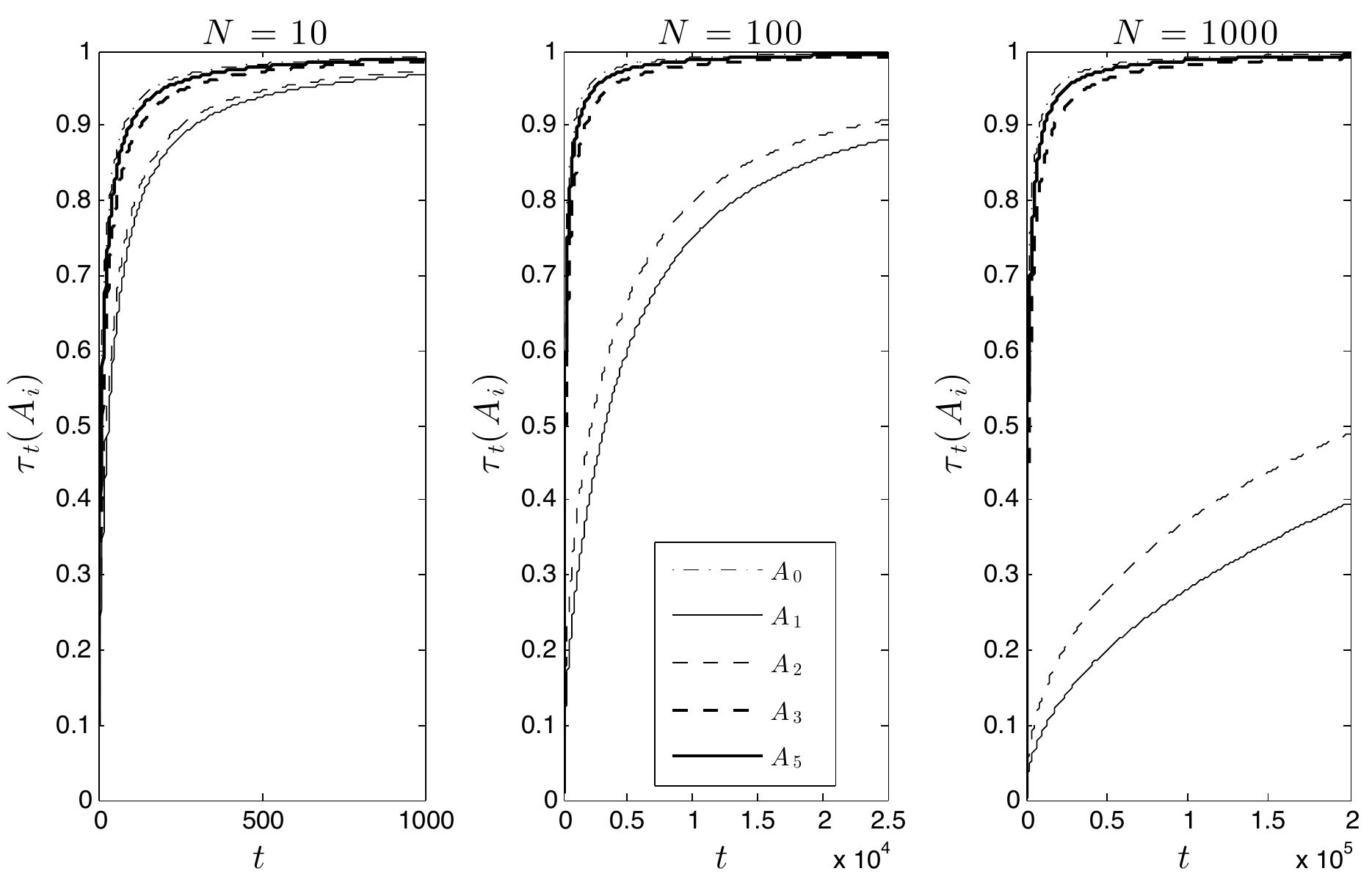}
\end{center}
\caption{Evolution of $\tau_t(A_i)$ with $t$ for different values of $N$, for $A = A_0$, $A_1$, $A_2$ as before with the addition of $3$- and $5$-regular Ramanujan-type matrices $A_3$ and $A_5$.}\label{rama}
\end{figure}

The statistical efficiency of these graphs is further highlighted in Figure~\ref{rama}. It shows results for 3- and 5-regular Ramanujan-type matrices ($A_3$ and $A_5$) as well as the previous results for non-Ramanujan-type matrices $A_0$, $A_1$ and $A_2$ (see Figure~\ref{fig1}). 

We see that $A_3$ is already close to the statistical performance of $A_0$, the saturated network, and for all intents and purposes $A_5$ is essentially as good as $A_0$, even when there are $N=1000$ nodes; i.e., the statistical performance of the $5$-regular Ramanujan graph is barely distinguishable from that of the totally connected graph! Nevertheless, we must not forget that the possibility of building such efficient networks in real-world situations will ultimately depend on the specific application, and may not always be possible.

Next, assuming that the Ramanujan-type matrix $A$ is irreducible and aperiodic, it is apparent that there is a compromise to be made between the communication complexity of
the algorithm (as measured by the degree index $\mathscr C(A)=d$) 
and its statistical performance (as measured by the coefficient $\mathscr S(A)$).
Clearly, the two  are in conflict. Upon this a question arises: is it possible to reach a compromise in the range of statistical performances $\mathscr S(A)$ while varying the communication complexity between $d=3$ and $d=N$? 
The answer is affirmative, as shown in the following simulation exercise.

\begin{figure}[h]
\begin{center}
\includegraphics[width=12cm]{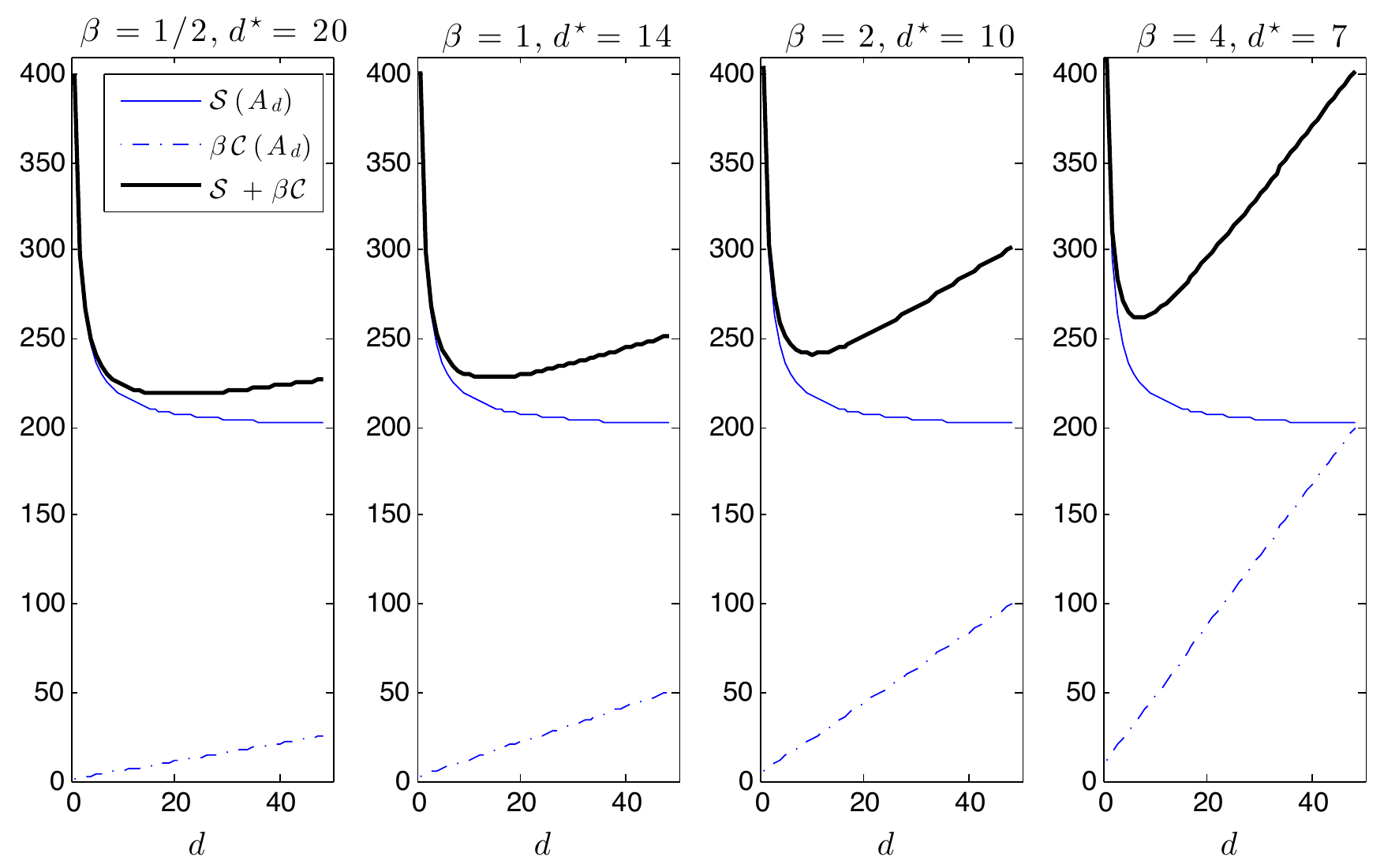}
\end{center}
\caption{Statistical efficiency vs communication complexity tradeoff for four different node communication penalties $\beta$. $d^{\star}$ is the $d$ which minimizes
$\mathscr S(A_d) + \beta \mathscr C(A_d)$.}\label{fig6}
\end{figure}

We fix $N=200$ and then for each $d = 3,\ldots,N$:
\begin{enumerate}[$(i)$]

\item Generate a matrix $A_d$ associated with a $d$-regular Ramanujan graph as before.

\item Compute the (non-unitary) eigenvalues $\gamma_2^{(d)}, \hdots, \gamma_N^{(d)}$ of the matrix $A_d$ and evaluate the sum $$\mathscr S(A_d)=\sum_{\ell=2}^N\frac{1}{1-\big(\gamma_{\ell}^{(d)}\big)^2} .$$

\item Plot $\mathscr S(A_d)$ and $\beta\mathscr C(A_d)=\beta d$ 
as well as penalized sums $\mathscr S(A_d) + \beta \mathscr C(A_d)$ for $\beta \in \{1/2, 1, 2, 4 \}$, where $\beta$ represents an explicit cost incurred when increasing the number of connections between nodes. 

\end{enumerate}
Results are shown in Figure~\ref{fig6}, where $d^{\star}$ refers to the $d$ for which the penalized sum $\mathscr S(A_d) + \beta \mathscr C(A_d)$ is minimized.
We observe that $\mathscr S(A_d)$ is decreasing whereas $\mathscr C(A_d)$ increases linearly.
The tradeoff between statistical efficiency and communication complexity can be seen as minimizing their penalized sum, where $\beta$ for example represents a monetary cost incurred by adding new network connections between nodes. We see that the optimal $d^{\star}$ and thus the number of node connections decreases as the cost of adding new ones increases.

Next, let us investigate the tradeoffs involved in the case where we have a large but fixed total number $T$ of data to be streamed to $N$ nodes, each receiving one new data value from time $t = 1$ to time $t = T/N$. In this context, the natural question to ask is how many nodes should we choose, and how much communication should we allow between them in order to get ``very good'' results for a ``low'' cost? Here a low cost comes from both limiting the number of nodes as well as the number of connections between them.

\begin{figure}[ht!]
\begin{center}
\includegraphics[width=8cm]{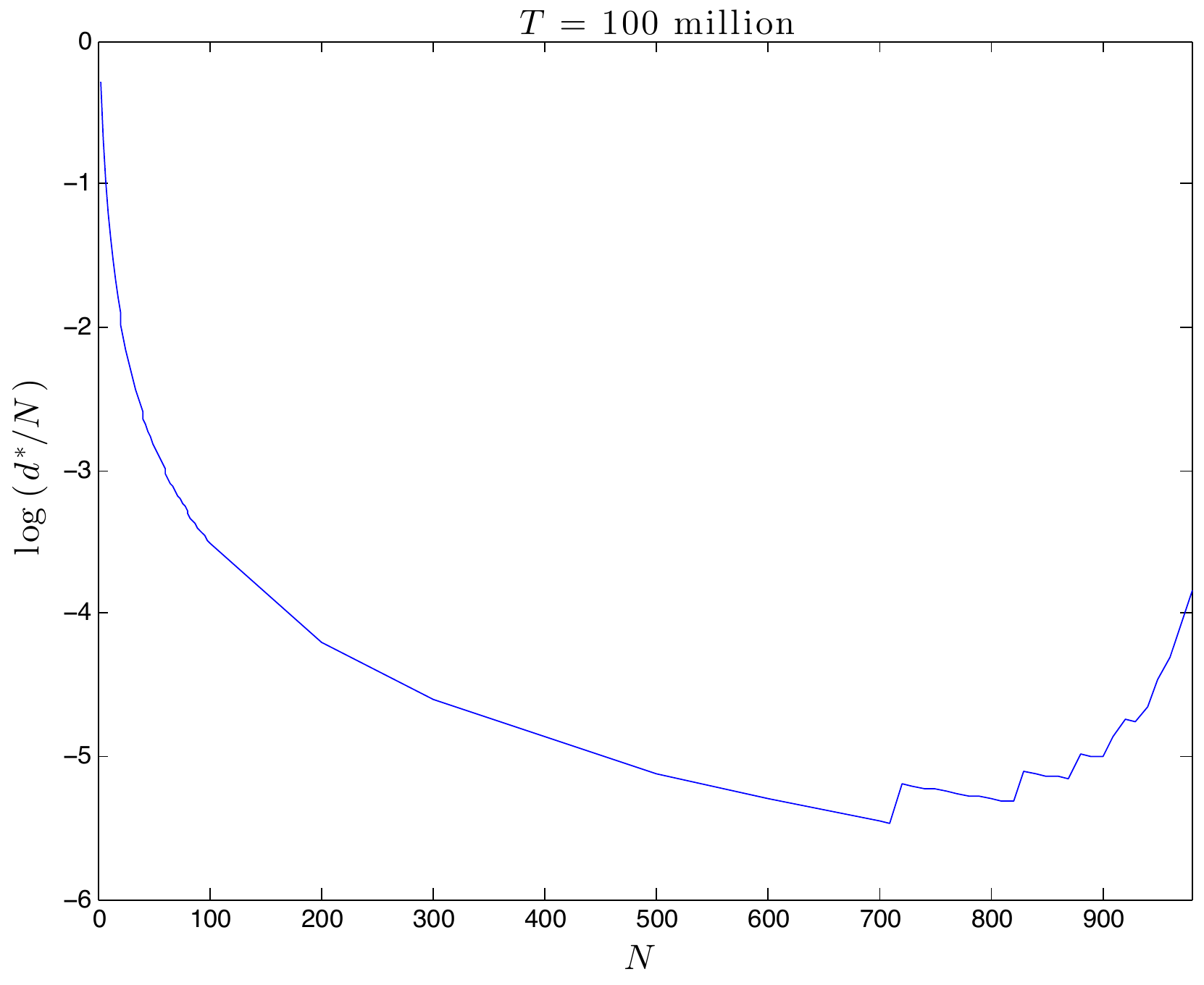}
\end{center}
\caption{Minimizing the number of nodes $N$ and the level of communication $d$ required between nodes to obtain a  performance ratio $\tau_t(A_{d}) \geq 0.99$ given a large fixed quantity of data $T$.}\label{fig3}
\end{figure}

In the same set-up for $A_d$ defined above, one way to look at this is to ask, for each $N$, what is the smallest $d \in \{3,\ldots,N\}$  and therefore the smallest communication cost $\mathscr C(A_d)=d$  for which the performance ratio $\tau_t(A_{d})$ is at least $0.99$ after receiving all the data, i.e., when $t = T/N$? Then, as there is also a cost associated with increasing $N$, minimizing  $\mathscr C(A_{d^{\star}})/N$ (where $d^{\star}$ is this smallest $d$ chosen) should help us choose the number of nodes $N$ and the amount of connection $\mathscr C(A_{d^{\star}})$ between them. The result of this is shown in Figure~\ref{fig3} for $T = 100$ million data points.
The minimum is found at $(N,d^{\star}) = (710,3)$, suggesting that with 100 million data points, one can get excellent performance results ($\tau_t(A_{d^{\star}}) \geq 0.99$) for a low cost with around 700 nodes, each connected only to three other nodes! Increasing $N$ further raises the cost necessary to obtain the same performance, both due to the price of adding more nodes, as well as requiring more connections between them: $d^{\star}$ must increase to 4, 5, and so on.
 
\section{Asynchronous models} 
The models considered so far assume that messages from one agent to another are immediately delivered. However, a distributed environment may be subject to communication delays, for instance when some processors compute faster than others or when latency and finite bandwidth issues perturb  message transmission. In the presence of such communication delays, it is conceivable that an agent will end up averaging its own value with an outdated value from another processor. Situations of this type fall within the framework of distributed asynchronous computation \citep[][]{TsBeAt86,BeTs89}. In the present section, we have in mind a model where  agents do not have to wait at predetermined moments for predetermined messages to become available. We thus allow some agents to compute faster and execute more iterations than others and allow communication delays to be substantial.

Communication delays are incorporated into our model as follows. For $B$ a nonnegative integer, we assume that the last instant before $t$ where agent $j$ sent a message to agent $i$ is $t-B_{ij}$, where $B_{ij}\in \{0,\hdots,B\}$. Put differently, recalling that $\hat \theta_t^{(i)}$ is the estimate held by agent $i$ at time $t$, we have
\begin{equation}
\label{retard}
\hat \theta_{t+1}^{(i)}=\frac{1}{t+1} \sum_{j=1}^N a_{ij} (t-B_{ij})\hat \theta^{(j)}_{t-B_{ij}}+\frac{1}{t+1}X_{t+1}^{(i)}, \quad t\geq 1.
\end{equation}
Thus, at time $t$, when agent $i$ uses the value of another agent $j$, this value is not necessarily the most recent one $\hat \theta_t^{(j)}$, but rather an outdated one $ \hat \theta_{t-B_{ij}}^{(j)}$, where $B_{ij}$ represents the communication delay. The time instants $t-B_{ij}$ are deterministic and, in any case, $0\leq B_{ij}\leq B$, i.e., we assume that  delays are bounded. Notice that some of the values $t-B_{ij}$ in (\ref{retard}) may be negative---in this case, by convention we set $\hat \theta^{(j)}_{t-B_{ij}}=0$. Our goal is to establish a counterpart to Theorem \ref{convergence} in the presence of communication delays. As usual, we set $\btheta_t=(\hat \theta_t^{(1)}, \hdots, \hat \theta^{(N)}_t)^\top$. 

Let $\kappa(t)$ be the smallest $\ell$ such that for all $(k_0,\hdots,k_\ell)\in \{1,\hdots,N\}^{\ell+1}$ satisfying 
$\prod_{j=1}^\ell a_{k_{j-1}k_j}>0,$
we have 
$$t-\ell-\sum_{j=1}^\ell B_{k_{j-1}k_j}\le B.$$
Observe that $t-\ell -\sum_{j=1}^\ell B_{k_{j-1}k_j}$ is the last time before $t$ when a message was sent from agent $k_0$ to agent $k_\ell$ via $k_1,\hdots,k_{\ell-1}$. Accordingly, $\kappa(t)$ is nothing but the smallest number of transitions needed to return at a time instant earlier than $B$, whatever the path. We note that $\kappa(t)$ is roughly of order $t$, since
$$\frac{1}{B+1}\le \liminf_{t\to\infty} \frac{\kappa(t)}{t}\le \limsup_{t\to\infty} \frac{\kappa(t)}{t}\le 1.$$

From now on, it is assumed that $A=A_1$, i.e., the irreducible, aperiodic, and symmetric matrix defined in (\ref{A1}). Besides its simplicity, this choice is motivated by the fact that $A_1$ is communication-efficient while its associated performance obeys
$$\tau_t(A)\approx 1-\frac{N^2}{6t}$$
for large $t$ and $N$. The main result of the section now follows.
\begin{theo}
\label{ASYNCHRO} 
Assume that $X$ is bounded and let $A=A_1$ be defined as in (\ref{A1}). Then, as $t \to \infty$,
$$\mathbb E \bigg\|\frac{t}{\kappa(t)} \btheta_t-\theta\mathbf 1\bigg\|^2=\emph{O}\bigg(\frac{1}{t}\bigg).$$
\end{theo}

The advantages one hopes to gain from asynchronism are twofold. First, a reduction of the synchronization penalty and a potential speed advantage over synchronous algorithms, perhaps at the expense of higher communication complexity. Second, a greater implementation flexibility and tolerance to system failure and uncertainty. On the other hand, the powerful result of Theorem \ref{ASYNCHRO} comes at the price of assumptions on the transmission network, which essentially demand that  communication delays $B_{ij}$ are time-independent. In fact, we find that the introduction of delays considerably complicates the consistency analysis of $\tau_t(A)$ even for the simple case of the empirical mean. This unexpected mathematical burden is due to the fact that the introduction of delays makes the analysis of the variance of the estimates quite complicated.
\section{Proofs}\label{proofs}
We start this section by recalling the following important theorem, whose proof can be found for example
in \citet[][Theorems 6.8.3 and 6.8.4]{FoFu04}. Here and elsewhere, $A$ stands for the stochastic
communication matrix.
\begin{theo}
\label{FF}
Let $\lambda_1,\hdots,\lambda_d$ be the eigenvalues of $A$ of unit modulus (with $\lambda_1=1$) and  $\Gamma$ be the set of eigenvalues of $A$ of modulus strictly smaller than 1. 
\begin{enumerate}[$(i)$]
\item There exist projectors $Q_1,\hdots,Q_d$ such that, for all $k\ge N$, 
$$A^k=\sum_{\ell=1}^d \lambda_{\ell}^k Q_{\ell}+\sum_{\gamma\in\Gamma} \gamma^k Q_\gamma(k),$$
where the matrices $\{Q_\gamma(k): k\ge N, \gamma\in\Gamma\}$ satisfy $Q_\gamma(k)Q_{\gamma'}(k')=Q_\gamma(k+k')$ if $\gamma=\gamma'$, and $0$ otherwise. In addition, for all $\gamma \in \Gamma$, $\lim_{k \to \infty} \gamma^kQ_{\gamma}(k)=0$.
\item The sequence $(A^k)_{k\geq 0}$ converges in the Ces\`aro sense to $Q_1$, i.e., 
$$ \frac{1}{t} \sum_{k=0}^t A^k \to Q_1 \quad \mbox{as $t\to\infty$}.$$
\end{enumerate}
\end{theo}
\subsection{Proof of Proposition \ref{PROPO1}}
According to (\ref{BLVDP}), since $A^k$ is a stochastic matrix, we have
$$\btheta_{t}-\theta \mathbf 1= \frac{1}{t}\sum_{k=0}^{t-1}A^k(\bX_{t-k}-\theta \mathbf 1).$$
Therefore, it may be assumed, without loss of generality, that $\theta=0$. Thus,
$$\tau_t(A)=\frac{\mathbb E \left\| \bar{\mathbb X}_{Nt}\mathbf 1\right\|^2}{\mathbb E\|\btheta_t \|^2}.$$
Next, let $A^k=(a_{ij}^{(k)})_{1\leq i,j\leq N}$. Then, for each $i \in\{1, \hdots, N\}$,
$$\hat \theta_t^{(i)}=\frac{1}{t}\sum_{k=0}^{t-1}\sum_{j=1}^N a_{ij}^{(k)} X_{t-k}^{(j)}, \quad t\geq 1.$$
By independence of the samples,
$$\mathbb E \big(\hat \theta_t^{(i)}\big)^2=\frac{\sigma^2}{t^2}\sum_{k=0}^{t-1} \sum_{j=1}^N \big(a_{ij}^{(k)}\big)^2.$$
Upon noting that $\mathbb E (\bar{\mathbb X}_{Nt})^2=\frac{\sigma^2}{Nt}$, we get
\begin{align*}
\tau_t(A) & =\frac{N \mathbb E \big(\bar{\mathbb X}_{Nt}\big)^2}{\mathbb E \big(\hat \theta_t^{(1)}\big)^2+\cdots+\mathbb E \big(\hat \theta_t^{(N)}\big)^2}\\
& =\frac{t}{\sum_{k=0}^{t-1} \|A^k\|^2}.
\end{align*}
Since each $A^k$ is a stochastic matrix, $\|A^k\|^2\le N$ and, by the Cauchy-Schwarz inequality, $\|A^k\|\ge 1$. Thus, 
$\frac{1}{N}\le \tau_t(A)\le 1$, the lower bound being achieved when $A$ is the identity matrix. 

Let us now assume that $A$ is reducible, and let $C\subsetneq\{1, \hdots , N\}$ be a recurrence class. Arguing as above, we obtain that for all $i\in C$, 
$$\mathbb E\big(\hat \theta^{(i)}_t\big)^2=\frac{\sigma^2}{t^2} \sum_{k=0}^{t-1} \sum_{j=1}^N \big(a^{(k)}_{ij}\big)^2\ge \frac{\sigma^2}{t^2} \sum_{k=0}^{t-1} \sum_{j\in C} \big(a^{(k)}_{ij}\big)^2.$$
Since $C$ is a recurrence class, the restriction of $A$ to entries in $C$ is a stochastic matrix as well. Thus, setting $N_1=|C|$, by the Cauchy-Schwarz inequality,
\begin{equation*} 
\mathbb E \big(\hat \theta^{(i)}_t\big)^2 \geq \left\{
\begin{array}{ll}
 \frac{\sigma^2}{tN_1} & \mbox{ if $i\in C$}\\
 \frac{\sigma^2}{tN}& \mbox{ otherwise.}
\end{array}
\right.
\end{equation*}
To conclude,
\begin{align*}
\tau_t(A) & =  \frac{\sigma^2/t}{\sum_{i\in C} \mathbb E \big(\hat \theta^{(i)}_t\big)^2+\sum_{i\notin C} \mathbb E \big(\hat \theta^{(i)}_t\big)^2}\\
& \le  \frac{1}{1+(N-N_1)/N}\\
&\le  \frac{N}{N+1},
\end{align*}
since $N-N_1\geq 1$.
\subsection{Proof of Lemma \ref{LEM1}}
As in the previous proof, we assume that $\theta=0$. Recall that
$$\btheta_{t}=\frac{1}{t} \sum_{k=0}^{t-1} A^k \bX_{t-k}, \quad t\geq 1.$$
Thus, for all $t\geq 1$, 
\begin{align*}
\mathbb E\|\btheta_t\|^2&=\frac{1}{t^2}\mathbb E \bigg \| \sum_{k=0}^{t-1} A^k \bX_{t-k}\bigg\|^2\nonumber\\
&=\frac{1}{t^2}\sum_{k=0}^{t-1} \mathbb E   \| A^k \bX_{t-k}\|^2\nonumber\\
& \quad \mbox{(by independence of $\bX_1, \hdots, \bX_t$)}\nonumber\\
&=\frac{1}{t^2}\mathbb E\bX_1^{\top} \bigg(\sum_{k=0}^{t-1}(A^k)^{\top}A^k\bigg)\bX_1.
\end{align*}
Denote by $\lambda_1=1, \hdots, \lambda_d$ the eigenvalues of $A$ of modulus 1, and let $\Gamma$ be the set of eigenvalues $\gamma$ of $A$ of modulus strictly smaller than 1. According to Theorem \ref{FF}, there exist projectors $Q_1,\hdots, Q_d$ and matrices $Q_{\gamma}(k)$ such that for all $k\geq N$,
$$A^k=\sum_{\ell=1}^d \lambda_{\ell}^k Q_{\ell}+\sum_{\gamma \in \Gamma} \gamma^k Q_{\gamma}(k) .$$
Therefore,
\begin{align*}
\sum_{k=0}^{t-1} (A^k)^{\top}A^k &=\sum_{k=0}^{t-1} (\bar A^k)^{\top}A^k\\
&=\sum_{k=0}^{t-1}\bigg( \sum_{\ell=1}^d \bar \lambda_{\ell}^k \bar Q_{\ell}+\sum_{\gamma \in \Gamma}\bar \gamma^k \bar Q_{\gamma}(k)\bigg)^{\top}\bigg( \sum_{j=1}^d \lambda_{j}^k Q_{j}+\sum_{\gamma \in \Gamma}\gamma^k Q_{\gamma}(k)\bigg)\\
& =\sum_{k=0}^{t-1}\sum_{\ell, j=1}^d \bar \lambda_{\ell}^k\lambda_{j}^k \bar Q_{\ell}^{\top} Q_{j}+\mbox{o}(t).
\end{align*}
Here, we have used Ces\`aro's lemma combined with the fact that for any $\gamma \in \Gamma$, $\lim_{k \to \infty}\gamma^k Q_{\gamma}(k)= 0$ (Theorem \ref{FF}). 

Since $A$ is irreducible, according to the Perron-Frobenius theorem \citep[e.g.,][page 240]{GrSt01}, we have that $\lambda_{\ell}=e^{\frac{2\pi i(\ell-1)}{d}}$, $1\leq \ell \leq d$. Accordingly,
$$\bar \lambda_{\ell}\lambda_j=e^{\frac{2\pi i (j-\ell)}{d}}=1 \Leftrightarrow j=\ell.$$ 
Thus,
$$\sum_{k=0}^{t-1}(A^k)^{\top} A^k=t\sum_{\ell=1}^d \bar Q_{\ell}^{\top} Q_{\ell}+\mbox{O}(1)+\mbox{o}(t).$$
Letting $Q=\sum_{\ell=1}^d \bar Q_{\ell}^{\top} Q_{\ell}$, we obtain
\begin{align}
t\mathbb E \| \btheta _t\|^2 & =\mathbb E \bX_1^{\top}Q \bX_1+\mathbb E \bX_1^{\top}\bigg(\frac{1}{t}\sum_{k=0}^{t-1} (A^k)^{\top} A^k-Q\bigg)\bX_1\label{CNU26}\\
& =\mathbb E \bX_1^{\top} Q \bX_1+\mbox{o}(1) \nonumber\\
&=\sum_{\ell=1}^d \mathbb E \| Q_{\ell}\bX_1\|^2+\mbox{o}(1). \nonumber
\end{align}
Denoting by $Q_{\ell,ij}$ the $(i,j)$-entry of $Q_{\ell}$, we conclude
\begin{align*}
t\mathbb E \| \btheta _t\|^2& =\sum_{\ell=1}^d \mathbb E \sum_{i=1}^N\bigg(\sum_{j=1}^N Q_{\ell,ij}X_1^{(j)}\bigg)^2+\mbox{o}(1)\\
&=\sigma^2\sum_{\ell=1}^d  \sum_{i,j=1}^N Q^2_{\ell,ij}+\mbox{o}(1)\\
& \quad \mbox{(by independence of $X_1^{(1)}, \hdots, X_1^{(N)}$)}\\
&=\sigma^2\sum_{\ell=1}^d\|Q_{\ell}\|^2+\mbox{o}(1).
\end{align*}
Lastly, recalling that $\mathbb E \| \bar{\mathbb X}_{Nt}\mathbf 1\|^2=\frac{\sigma^2}{t}$, we obtain
$$\tau_t(A)=\frac{1}{\sum_{\ell=1}^d\|Q_{\ell}\|^2+\mbox{o}(1)}=\frac{1}{\sum_{\ell=1}^d\|Q_{\ell}\|^2}+\mbox{o}(1).$$
\subsection{Proof of Theorem \ref{convergence}}

{\bf Sufficiency.} Assume that $A$ is irreducible, aperiodic, and bistochastic. The first two conditions imply that $1$ is the unique eigenvalue of $A$ of unit modulus. Therefore, according to Lemma \ref{LEM1}, we only need to prove that the projector $Q_1$ satisfies $\|Q_1\|=1$. 

Since $A$ is bistochastic, its stationary distribution is the uniform distribution on $\{1,\hdots,N\}$. Moreover, since $A$ is irreducible and aperiodic, we have, as $k\to\infty$,
$$A^k \to\frac{1}{N} \begin{pmatrix} 1 & 1 &\hdots & 1\\ \vdots & \vdots & \vdots &\vdots \\ 1 & 1 &\hdots & 1\end{pmatrix}.$$
By comparing this limit with that of the second statement of Theorem \ref{FF}, we conclude by Ces\`aro's lemma that 
$$
Q_1= \frac{1}{N}\begin{pmatrix} 1 & 1 & \hdots& 1\\ \vdots & \vdots & \vdots &\vdots \\ 1 & 1 & \hdots & 1\end{pmatrix}.
 $$
This implies in particular that $\|Q_1\|=1$. 

{\bf Necessity}. Assume that $\tau_t(A)$ tends to 1 as $t\to\infty$. According to Proposition \ref{PROPO1}, $A$ is irreducible. Thus, by Lemma \ref{LEM1}, we have $\sum_{\ell=1}^d \|Q_\ell\|^2=1$. Observe, since each $Q_{\ell}$ is a projector, that $\|Q_{\ell}\|\geq 1$. Therefore, the identity $\sum_{\ell=1}^d \|Q_{\ell}\|^2=1$ implies $d=1$ and $\|Q_1\|=1$. We conclude that $A$ is aperiodic. 

Then, since $A$ is irreducible and aperiodic, we have, as $k\to \infty$,
$$A^k \to \begin{pmatrix}\boldsymbol{\mu}\\ \vdots \\ \boldsymbol{\mu}\end{pmatrix},$$
where $\boldsymbol{\mu}$ is the stationary distribution of $A$, represented as a row vector. Comparing once again this limit with the second statement of Theorem \ref{FF}, we see that 
$$
Q_1= \begin{pmatrix}\boldsymbol{\mu}\\ \vdots \\ \boldsymbol{\mu}\end{pmatrix}.
 $$
Thus, $\|Q_1\|^2=N\| {\boldsymbol{\mu}}\|^2=1$. In particular, letting $\boldsymbol{\mu}=(\mu_1, \hdots, \mu_N)$, we have 
$$N\sum_{i=1}^N \mu_i^2=\sum_{i=1}^N \mu_i.$$
This is an equality case in the Cauchy-Schwarz inequality, from which we deduce that $\boldsymbol{\mu}$ is the uniform distribution on $\{1,\hdots,N\}$. Since $\boldsymbol{\mu}$ is  the stationary distribution of $A$, this implies that $A$ is bistochastic.

\subsection{Proof of Proposition \ref{CORO1}}
If $A$ is irreducible and aperiodic, then by Lemma \ref{LEM1}, $\tau_t(A)\to \frac{1}{\|Q_1\|^2}$ as $t\to \infty$. But, as $k \to \infty$,
$$A^k \to \begin{pmatrix}\boldsymbol{\mu}\\ \vdots \\ \boldsymbol{\mu}\end{pmatrix},$$
where the stationary distribution $\boldsymbol{\mu}$ of $A$ is represented as a row vector. By the second statement of Theorem \ref{FF}, we conclude that $\|Q_1\|^2=N\|\boldsymbol{\mu}\|^2$.

\subsection{Proof of Theorem \ref{theovitesse}}
Without loss of generality,  assume that $\theta=0$. Since $A$ is irreducible and aperiodic, the matrix $Q$ in the proof of Lemma \ref{LEM1} is $Q=Q_1^\top Q_1$. Moreover, since $A$ is also bistochastic, we have already seen that as $k\to\infty$,
\begin{equation}
\label{++}
A^k\to \frac{1}{N}\begin{pmatrix} 1 & 1 &\hdots & 1\\ \vdots & \vdots & \vdots &\vdots \\ 1 & 1 & \hdots& 1\end{pmatrix}.
\end{equation}
However, by the second statement of Theorem \ref{FF}, the above matrix is equal to $Q_1$. Thus, the projector $Q_1$ is symmetric, which implies $Q=Q_1$. 

Next, we deduce from (\ref{CNU26}) that
\begin{align}
\tau_t(A) & =  \frac{\sigma^2}{\mathbb E \bX_1^\top Q \bX_1+\mathbb E \bX_1^\top \big(\frac{1}{t}\sum_{k=0}^{t-1} (A^k)^\top A^k-Q\big)\bX_1}\nonumber\\
& =  \frac{\sigma^2}{\sigma^2+\mathbb E \bX_1^{\top} \big(\frac{1}{t}\sum_{k=0}^{t-1} A^{2k}-Q\big)\bX_1},\label{bla2}
\end{align}
by symmetry of $A$ and the fact that $\mathbb E \bX_1^\top Q \bX_1=\sigma^2$. The symmetric matrix $A$ can be put into the form
$$A=UDU^{\top},$$
where $U$ is a unitary matrix with real entries (so, $U^{\top}=U^{-1}$) and
$D=\mbox{diag}(1,\gamma_2, \hdots, \gamma_N)$, with $1> \gamma_2 \geq \cdots\geq \gamma_N>-1$. Therefore, as $k\to \infty$,
$$\frac{1}{t}\sum_{k=0}^{t-1} A^{2k}=U \bigg( \frac{1}{t}\sum_{k=0}^{t-1} D^{2k}\bigg)U^{\top}\to U \begin{pmatrix} 1 & 0 & \hdots & 0\\ 0 &  0 &\hdots & 0\\ \vdots & \vdots & \vdots & \vdots\\0 & 0 & \hdots & 0\end{pmatrix}U^\top.$$
However, by (\ref{++}) and Ces\`aro's lemma,
$$\frac{1}{t}\sum_{k=0}^{t-1} A^{2k}\to  Q\quad \mbox{as $k\to \infty$}.$$
It follows that $Q=UMU^{\top}$,
where
$$M=
\begin{pmatrix} 1 & 0 & \hdots & 0\\ 0 &  0 &\hdots  & 0\\ \vdots & \vdots & \vdots & \vdots\\0 & 0 & \hdots& 0\end{pmatrix}.$$
Thus,
\begin{align*}
\frac{1}{t}\sum_{k=0}^{t-1}A^{2k} -Q&=U \bigg( \frac{1}{t} \sum_{k=0}^{t-1} D^{2k}-M\bigg)U^{\top}\\
&=U\bigg( \frac{1}{t}\sum_{k=0}^{t-1}\mbox{diag}\big(0,\gamma_2^{2k}, \hdots, \gamma_N^{2k}\big)\bigg)U^{\top}\\
&=U\mbox{diag}\bigg(0,\frac{1}{t}\frac{1-\gamma_2^{2t}}{1-\gamma_2^2}, \hdots, \frac{1}{t}\frac{1-\gamma_N^{2k}}{1-\gamma_N^2}\bigg)U^{\top}.
\end{align*}

Next, set
$$\alpha_{\ell}=\frac{1}{t}\frac{1-\gamma_{\ell}^{2t}}{1-\gamma_{\ell}^2}, \quad 2 \leq \ell \leq N,$$
and let $U=(u_{ij})_{1\leq i,j\leq N}$. With this notation, the $(i,j)$-entry of the matrix $\frac{1}{t}\sum_{k=0}^{t-1}A^{2k} -Q$   is
$$\sum_{\ell=2}^N u_{i\ell} \alpha_{\ell}u_{j\ell}.$$
Hence,
$$
\bX_1^{\top}\bigg(\frac{1}{t}\sum_{k=0}^{t-1}A^{2k} -Q\bigg)\bX_1 =\sum_{i=1}^N X_1^{(i)}\sum_{j=1}^N\bigg(\sum_{\ell=2}^N u_{i \ell} \alpha_{\ell} u_{j \ell}\bigg)X_1^{(j)}.$$
Thus,
\begin{align*}
\mathbb E\bX_1^{\top}\bigg(\frac{1}{t}\sum_{k=0}^{t-1}A^{2k} -Q\bigg)\bX_1 &=\sigma^2 \sum_{i=1}^N \sum_{\ell=2}^N u_{i\ell}\alpha_{\ell}u_{i\ell}\\
&=\sigma^2 \sum_{i=1}^N \sum_{\ell=2}^N \alpha_{\ell}u_{i\ell}^2\\
&=\sigma^2 \sum_{\ell=2}^N \alpha_{\ell}\\
&=\frac{\sigma^2}{t}\sum_{\ell=2}^N\frac{1-\gamma_{\ell}^{2t}}{1-\gamma_{\ell}^2}.
\end{align*}
We conclude from (\ref{bla2}) that 
$$
\tau_t(A) = \frac{1}{1 + \frac{1}{t}\sum_{\ell = 2}^{N}\frac{1 - \gamma_{\ell}^{2t}}{1 - \gamma_{\ell}^2}}.
$$
This shows the first statement of the theorem. Using the inequality $\frac{1}{1+x}\ge 1-x$, valid for all $x\ge 0$, we have
\begin{align*}
\tau_t(A) & \ge  1-\frac{1}{t} \sum_{\ell=2}^N \frac{1-\gamma_\ell^{2t}}{1-\gamma_\ell^2}\\
& \ge  1-\frac{\mathscr S(A)}{t}.
\end{align*}
Finally, evoking the inequality $\frac{1}{1+x}\le 1-x+x^2$, valid for all $x\ge 0$, we conclude
\begin{align*}
\tau_t(A) & \le  1-\frac{1}{t} \sum_{\ell=2}^N \frac{1-\gamma_\ell^{2t}}{1-\gamma_\ell^2}+\bigg(\frac{1}{t} \sum_{\ell=2}^N \frac{1-\gamma_\ell^{2t}}{1-\gamma_\ell^2}\bigg)^2\\
& \le 1 -\frac{\mathscr S(A)}{t}+\Gamma^{2t}(A) \frac{\mathscr S(A)}{t}+\Big(\frac{\mathscr S(A)}{t}\Big)^2.
\end{align*}
\subsection{Proof of Theorem \ref{ASYNCHRO}} 
From now on, we fix $k_0\in\{1,\hdots,N\}$ and let $Z_t^{(i)}=t\hat \theta_t^{(i)}$ for any $i\in\{1,\hdots,N\}$. Thus, for all $t \geq 1$,
$$Z_t^{(k_0)}=\sum_{k=1}^N a_{k_0k}Z_{t-B_{k_0k}-1}^{(k)}+X_t^{(k_0)},$$
and
\begin{equation}
\label{ordre2}
Z_t^{(k_0)} = \sum_{k_1,k_2=1}^N a_{k_0k_1}a_{k_1k_2}Z_{t-B_{k_0k_1}-B_{k_1k_2}-2}^{(k_2)} + \sum_{k_1=1}^N a_{k_0k_1} X_{t-B_{k_0k_1}-1}^{(k_1)}+X_t^{(k_0)}.
\end{equation}
Our first task is to iterate this formula. To do so, we need  additional notation. For $\ell$ a positive integer and $k\in\{1,\hdots,N\}$,  let $\underline K^\ell(k)$ be the set of vectors in $\{1,\hdots,N\}^{\ell+1}$ of the form $(k_0,k_1,\hdots,k_{\ell-1},k)$ such that $w(\underline K^\ell(k))>0$, where
$$w\big(\underline K^\ell(k)\big)=a_{k_0k_1}a_{k_1k_2}\hdots a_{k_{\ell-2}k_{\ell-1}}a_{k_{\ell-1} k}.$$
In particular, by our choice of $A$, we have $w(\underline K^\ell(k))=2^{-\ell}$ for any $k$. Next, we set $$\Delta\big(\underline K^\ell(k)\big)=\ell+B_{k_0k_1}+B_{k_1k_2}+\cdots+B_{k_{\ell-2}k_{\ell-1}}+B_{k_{\ell-1} k}.$$
When  $\ell=0$, then by convention $\underline K^0(k)=(k_0)$, $w(\underline K^0(k))=1$ if $k=k_0$ and $0$ otherwise, and $\Delta(\underline K^0(k))=0$. 

We are now ready to iterate (\ref{ordre2}). To do so, observe that
\begin{align}
Z_t^{(k_0)} &=  \sum_{k=1}^N \sum_{\underline K^{\kappa(t)}(k)} 
w\big(\underline K^{\kappa(t)}(k)\big) Z^{(k)}_{t-\Delta(\underline K^{\kappa(t)}(k))}\nonumber \\ 
&  \quad + \sum_{\ell=0}^{\kappa(t)-1} \sum_{k=1}^N  \sum_{\underline K^{\ell}(k)} w\big(\underline K^{\ell}(k)\big) X^{(k)}_{t-\Delta(\underline K^\ell (k))}\nonumber\\
& \stackrel{\mbox{\tiny def}}{=} R_t^1+R_t^2.\label{horreur}
\end{align}
By the definition of $\kappa(t)$, for all $k\in\{1,\hdots,N\}$, $t-\Delta(\underline K^{\kappa(t)}(k))\le B$. Since $X$ is bounded, we deduce that there exists $C>0$ such that 
$$|R_t^1|\le C \sum_{k=1}^N \sum_{\underline K^{\kappa(t)}(k)} 
w\big(\underline K^{\kappa(t)}(k)\big).$$
This implies that $|R_t^1|\le C$. To see this, note that $A^{\kappa(t)}$ is a stochastic matrix and that for all $k\in\{1,\hdots,N\}$, 
$$\sum_{\underline K^{\kappa(t)}(k)} w\big(\underline K^{\kappa(t)}(k)\big)=(A^{\kappa(t)})_{k_0k}.$$
The analysis of the term $R_t^2$ is more delicate. The difficulty arises from the fact that this term is \emph{not} a sum of independent random variables, and therefore its components must be grouped. Since each $B_{ij}$ is smaller than $B$ and $\Delta(\underline K^\ell(k))=x$ implies $x\ge \ell$, we obtain 
\begin{align*}
R_t^2 & =  \sum_{\ell=0}^{\kappa(t)-1} \sum_{k=1}^N \sum_{x=0}^{(B+1)\ell} \sum_{\underline K^\ell(k) : \Delta(\underline K^\ell(k))=x}w\big(\underline K^\ell(k)\big)\, X_{t-x}^{(k)}\\
& =  \sum_{x=0}^{(B+1)(\kappa(t)-1)} \sum_{k=1}^N \sum_{\ell=\lfloor x/(B+1)\rfloor+1}^{x} \sum_{\underline K^\ell(k) : \Delta(\underline K^\ell(k))=x} w\big(\underline K^\ell(k)\big) \, X_{t-x}^{(k)}
\end{align*}
($\lfloor \cdot \rfloor$ is the floor function). By independence of the $X_j^{(i)}$, we get 
$${\rm Var}(R_t^2)=\sigma^2 \sum_{x=0}^{(B+1)(\kappa(t)-1)} \sum_{k=1}^N \bigg(\sum_{\ell=\lfloor x/(B+1) \rfloor +1}^{x} \sum_{\underline K^\ell(k) : \Delta(\underline K^\ell(k))=x} w\big(\underline K^\ell(k)\big)\bigg)^2.$$
Recalling that $w(\underline K^\ell(k))=2^{-\ell}$, we obtain
$${\rm Var}(R_t^2)=\sigma^2 \sum_{x=0}^{(B+1)(\kappa(t)-1)} \sum_{k=1}^N \bigg(\sum_{\ell= \lfloor x/(B+1) \rfloor+1}^{x} \frac{1}{2^\ell}\,  \Big| \underline K^\ell(k) : \Delta\big(\underline K^\ell(k)\big)=x\Big| \bigg)^2.$$
Next, consider the Markov chain $(Y_n)_{n\ge 0}$ with transition matrix $A$ such that $Y_0=k_0$. Observe that 
$$\mathbb P\Big(Y_\ell=k,\sum_{j=1}^\ell B_{Y_{j-1}Y_j}=x-\ell\Big)=\frac{1}{2^\ell}\, \Big| \underline K^\ell(k) :\Delta\big(\underline K^\ell(k)\big)=x\Big |.$$
Moreover, for fixed $x$, the events 
$$\bigg\{\sum_{j=1}^\ell B_{Y_{j-1}Y_j}=x-\ell\bigg\}, \quad \Big \lfloor \frac{x}{B+1}\Big \rfloor+1 \leq \ell \leq  x,$$
are disjoint since the $B_{ij}$ are nonnegative. Thus, 
$$\sum_{\ell=\lfloor x/(B+1) \rfloor+1}^{x} \frac{1}{2^\ell} \,\Big| \underline K^\ell(k) : \Delta\big(\underline K^\ell(k)\big)=x\Big|\le 1,$$
and so,
\begin{equation}
\label{bla1} {\rm Var}(R_t^2)\le \sigma^2 \sum_{x=0}^{(B+1)(\kappa(t)-1)} \sum_{k=1}^N 1=\sigma^2N\big((B+1)\kappa(t)-B\big).
\end{equation}
The expectation of $R_t^2$ is easier to compute. Indeed, since each $A^\ell$ is a stochastic matrix,
$$
\mathbb E R_t^2 =\theta \sum_{\ell=0}^{\kappa(t)-1} \sum_{k=1}^N \sum_{\underline K^\ell(k)} w\big(\underline K^\ell(k)\big) = 
\theta \sum_{\ell=0}^{\kappa(t)-1} \sum_{k=1}^N (A^{\ell})_{k_0k}= \theta \kappa(t) .$$
Combining (\ref{horreur}), (\ref{bla1}), and the fact that $|R_t^1|\le C$, we obtain 
\begin{align*}
\mathbb E \bigg(\frac{t}{\kappa(t)} \hat \theta_t^{(k_0)}-\theta\bigg)^2 & =  \mathbb E \bigg(\frac{R_t^1}{\kappa(t)}+\frac{R_t^2}{\kappa(t)}-\theta\bigg)^2\\
& =  \mathbb E \bigg( \frac{R_t^2-\mathbb E R_t^2}{\kappa(t)}+\frac{R_t^1}{\kappa(t)}\bigg)^2\\
& = \mbox{O}\bigg(\frac{1}{\kappa(t)}\bigg).
\end{align*}
The result follows from the identity $1/\kappa(t)=\mbox{O}(1/t)$. 
\bibliography{biblio-bbc}

\begin{thebibliography}{26}
\providecommand{\natexlab}[1]{#1}
\providecommand{\url}[1]{\texttt{#1}}
\expandafter\ifx\csname urlstyle\endcsname\relax
  \providecommand{\doi}[1]{doi: #1}\else
  \providecommand{\doi}{doi: \begingroup \urlstyle{rm}\Url}\fi

\bibitem[Alon(1986)]{Al86}
N.~Alon.
\newblock Eigenvalues and expanders.
\newblock \emph{Combinatorica}, 6:\penalty0 83--96, 1986.

\bibitem[Bertsekas and Tsitsiklis(1997)]{BeTs89}
D.P. Bertsekas and J.N. Tsitsiklis.
\newblock \emph{Parallel and Distributed Computation: Numerical Methods}.
\newblock Athena Scientific, Belmont, 1997.

\bibitem[Bianchi et~al.(2011{\natexlab{a}})Bianchi, Fort, Hachem, and
  Jakubowicz]{BiFoHaJa11}
P.~Bianchi, G.~Fort, W.~Hachem, and J.~Jakubowicz.
\newblock Convergence of a distributed parameter estimator for sensor networks
  with local averaging of the estimates.
\newblock In \emph{Proceedings of the 36th IEEE International Conference on
  Acoustics, Speech and Signal Processing}, 2011{\natexlab{a}}.

\bibitem[Bianchi et~al.(2011{\natexlab{b}})Bianchi, Fort, Hachem, and
  Jakubowicz]{BiFoHaJa11-2}
P.~Bianchi, G.~Fort, W.~Hachem, and J.~Jakubowicz.
\newblock Performance analysis of a distributed {R}obbins-{M}onro algorithm for
  sensor networks.
\newblock In \emph{Proceedings of the 19th European Signal Processing
  Conference}, 2011{\natexlab{b}}.

\bibitem[Bianchi et~al.(2013)Bianchi, Cl{\'e}men\c{c}on, Jakubowicz, and
  Morral]{BiClJaMo13}
P.~Bianchi, S.~Cl{\'e}men\c{c}on, J.~Jakubowicz, and G.~Morral.
\newblock On-line learning gossip algorithm in multi-agent systems with local
  decision rules.
\newblock In \emph{Proceedings of the 2013 IEEE International Conference on Big
  Data}, 2013.

\bibitem[Blondel et~al.(2005)Blondel, Hendrickx, Olshevsky, and
  Tsitsiklis]{BlHeOlTs05}
V.D. Blondel, J.M. Hendrickx, A.~Olshevsky, and J.N. Tsitsiklis.
\newblock Convergent in multiagent coordination, consensus, and flocking.
\newblock In \emph{Proceedings of the Joint 44th IEEE Conference on Decision
  and Control and European Control Conference}, 2005.

\bibitem[Bordenave(2015)]{Bo15}
C.~Bordenave.
\newblock A new proof of {F}riedman's second eigenvalue theorem and its
  extension to random lifts.
\newblock \emph{arXiv:1502.04482v1}, 2015.

\bibitem[Boyd et~al.(2006)Boyd, Ghosh, Prabhakar, and Shah]{BoGhPrSh06}
S.~Boyd, A.~Ghosh, B.~Prabhakar, and D.~Shah.
\newblock Randomized gossip algorithms.
\newblock \emph{IEEE Transactions on Information Theory}, 52:\penalty0
  2508--2530, 2006.

\bibitem[Duchi et~al.(2012)Duchi, Agarwal, and Wainwright]{DuAgWa12}
J.C. Duchi, A.~Agarwal, and M.J. Wainwright.
\newblock Dual averaging for distributed optimization: {C}onvergence analysis
  and network scaling.
\newblock \emph{IEEE Transactions on Automatic Control}, 57:\penalty0 592--606,
  2012.

\bibitem[Fiedler(1972)]{Fi72}
M.~Fiedler.
\newblock Bounds for eigenvalues of doubly stochastic matrices.
\newblock \emph{Linear Algebra and Its Applications}, 5:\penalty0 299--310,
  1972.

\bibitem[Foata and Fuchs(2004)]{FoFu04}
D.~Foata and A.~Fuchs.
\newblock \emph{Processus Stochastiques : {P}rocessus de {P}oisson,
  Cha\^{\i}nes de Markov et Martingales}.
\newblock Dunod, Paris, 2004.

\bibitem[Forero et~al.(2010)Forero, Cano, and Giannakis]{FoCaGi10}
P.A. Forero, A.~Cano, and G.B. Giannakis.
\newblock Consensus-based distributed support vector machines.
\newblock \emph{Journal of Machine Learning Research}, 11:\penalty0 1663--1707,
  2010.

\bibitem[Friedman(2008)]{Fr08}
J.~Friedman.
\newblock \emph{A Proof of AlonÕs Second Eigenvalue Conjecture and Related
  Problems}, volume 195 of \emph{Memoirs of the American Mathematical Society}.
\newblock American Mathematical Society, Providence, 2008.

\bibitem[Grimmett and Stirzaker(2001)]{GrSt01}
G.R. Grimmett and D.R. Stirzaker.
\newblock \emph{Probability and Random Processes. Third Edition}.
\newblock Oxford University Press, Oxford, 2001.

\bibitem[Jordan(2013)]{Jo13}
M.I. Jordan.
\newblock On statistics, computation and scalability.
\newblock \emph{Bernoulli}, 19:\penalty0 1378--1390, 2013.

\bibitem[Knight(2008)]{Kn08}
P.A. Knight.
\newblock The {S}inkhorn-{K}nopp algorithm: {C}onvergence and applications.
\newblock \emph{SIAM Journal on Matrix Analysis and Applications}, 30:\penalty0
  261--275, 2008.

\bibitem[Marcus et~al.(2015)Marcus, Spielman, and Srivastava]{MaSpSr15}
A.W. Marcus, D.A. Spielman, and N.~Srivastava.
\newblock Interlacing families {I}: {B}ipartite {R}amanujan graphs of all
  degrees.
\newblock \emph{Annals of Mathematics}, 182:\penalty0 307--325, 2015.

\bibitem[Mateos et~al.(2010)Mateos, Bazerques, and Giannakis]{MaBaGi10}
G.~Mateos, J.A. Bazerques, and G.B. Giannakis.
\newblock Distributed sparse linear regression.
\newblock \emph{IEEE Transactions on Signal Processing}, 58:\penalty0
  5262--5276, 2010.

\bibitem[Nilli(1991)]{Ni91}
A.~Nilli.
\newblock On the second eigenvalue of a graph.
\newblock \emph{Discrete Mathematics}, 91:\penalty0 207--210, 1991.

\bibitem[Olshevsky and Tsitsiklis(2011)]{OlTs11}
A.~Olshevsky and J.N. Tsitsiklis.
\newblock Convergence speed in distributed consensus and averaging.
\newblock \emph{SIAM Review}, 53:\penalty0 747--772, 2011.

\bibitem[Predd et~al.(2009)Predd, Kulkarni, and Poor]{PrKuPo09}
J.B. Predd, S.R. Kulkarni, and H.V. Poor.
\newblock A collaborative training algorithm for distributed learning.
\newblock \emph{{IEEE} Transactions on Automatic Control}, 55:\penalty0
  1856--1871, 2009.

\bibitem[Shlomo et~al.(2006)Shlomo, Linial, and Wigderson]{ShLiWi06}
H.~Shlomo, N.~Linial, and A.~Wigderson.
\newblock Expander graphs and their applications.
\newblock \emph{Bulletin of the American Mathematical Society}, 43:\penalty0
  439--561, 2006.

\bibitem[Sinkhorn and Knopp(1967)]{SiKn67}
R.~Sinkhorn and P.~Knopp.
\newblock Concerning nonnegative matrices and doubly stochastic matrices.
\newblock \emph{Pacific Journal of Mathematics}, 21:\penalty0 343--348, 1967.

\bibitem[Steger and Wormald(1999)]{St99}
A.~Steger and N.C. Wormald.
\newblock Generating random regular graphs quickly.
\newblock \emph{Combinatorics, Probability and Computing}, 8:\penalty0
  377--396, 1999.

\bibitem[Tsitsiklis et~al.(1986)Tsitsiklis, Bertsekas, and Athans]{TsBeAt86}
J.N. Tsitsiklis, D.P. Bertsekas, and M.~Athans.
\newblock Distributed asynchronous deterministic and stochastic gradient
  optimization algorithms.
\newblock \emph{{IEEE} Transactions on Automatic Control}, 31:\penalty0
  803--812, 1986.

\bibitem[Zhang et~al.(2013)Zhang, Duchi, and Wainwright]{ZhDuWa13}
Y.~Zhang, J.C. Duchi, and M.J. Wainwright.
\newblock Communication-efficient algorithms for statistical optimization.
\newblock \emph{Journal of Machine Learning Research}, 14:\penalty0 3321--3363,
  2013.

\end{thebibliography}

\end{document}